\documentclass[12pt,leqno]{amsart} 
\usepackage[boxed,ruled]{algorithm2e}
\newlength{\fixboxwidth}
\newcommand{\R}{{\mathbb R}}
\newcommand{\N}{{\mathbb N}}

\renewcommand{\rho}{{\varrho}}
\def\min{{\rm min}}

\def\a{{\alpha }} 
 
\def\e{{\varepsilon}} 
\def\alph{{\theta}}

\def\phi{{\varphi}}  % huebscher!? 
\def\om{\omega}
\newcommand{\vol}{\operatorname{vol}}
\newcommand{\rad}{\mathcal R^{\alpha}}
\newcommand{\radc}{\mathcal R_{C}}
\newcommand{\fad}{\mathcal F^\alpha}    %   neu 
\newcommand{\ball}{{B^d}}
\newcommand{\fo}{\mathcal F(\Omega)}
\newcommand{\fco}{\mathcal F_C(\Omega)}

\newcommand{\krd}{K_{\rho,\delta} }
\newcommand{\prd}{P_{\rho,\delta} }
\newcommand{\muo}{\mu_\Omega}
\newcommand{\mur}{\mu_\rho}
\newcommand{\dist}{\operatorname{dist}}
\newcommand{\vt}{S^{\mathrm{mean}}_n}
\newcommand{\vtm}{S^{\mathrm{mh}}_n}
\newcommand{\vtn}{S^{\mathrm{simple}}_n}
\newcommand{\lr}[1]{\left(#1\right)}
\newcommand{\abs}[1]{\left\vert #1 \right\vert} 
\newcommand{\norm}[2]{\Vert #1  \Vert _{#2}} 
\newcommand{\set}[1]{\left\{#1\right\}}
\newcommand{\expect}{\mathbf E}

\newcommand{\scalar}[2]{\langle #1,#2\rangle}

\theoremstyle{plain}
\newtheorem{theorem}{Theorem}
\newtheorem{lemma}{Lemma}
\newtheorem{proposition}{Proposition}
\newtheorem{corollary}{Corollary}

\theoremstyle{definition}
\newtheorem{rem}{Remark}
\begin{document}

\title{Simple Monte Carlo and the Metropolis algorithm} 
\author{Peter Math\'e} 
\address{Weierstrass Institute for Applied Analysis and 
Stochastics, Mohrenstrasse 39, D-10117 Berlin, Germany}
\email{mathe@wias-berlin.de}
\author{Erich Novak}
\address{ Friedrich Schiller University Jena, 
Mathem. Institute,
Ernst-Abbe-Platz 2, 
D-07743 Jena, Germany}
\email{novak@math.uni-jena.de}
\date{Version: \today}

\maketitle
\begin{center}
{\sl\large Dedicated to our dear colleague and friend Henryk
Wo\'zniakowski on the occasion of his 60th birthday. }  
\end{center}

\begin{abstract}
We study the integration of functions with 
respect to an unknown density.
% which is known only up to the normalizing factor. 
Information is available as
oracle calls to the integrand and to the non-normalized density 
function.
We are interested in analyzing the integration error
of optimal algorithms 
(or the complexity of the problem) with emphasis on
the variability of the weight function. 
For a corresponding large
class of problem instances we show that the complexity
grows linearly in the variability, and the simple Monte Carlo method
provides an almost optimal algorithm.
Under additional geometric restrictions (mainly log-concavity)
for the density
functions, we establish that a suitable adaptive
local Metropolis algorithm is almost optimal and 
outperforms any non-adaptive algorithm. 
\end{abstract}

\section{Introduction, Problem description}\label{s1} 
In many applications one wants to compute an integral of the form
\begin{equation}
\label{eq:base}
\int_\Omega f(x) \cdot c \rho(x) \, \mu(dx)  
\end{equation}
with a density $c \rho(x),\ x\in \Omega$, where $c >0$ is unknown
and $\mu$ is a probability measure. 
Of course we have 
$
{1}/{c} = \int_\Omega  \rho(x) \, \mu(dx),
$
but the numerical computation of the 
latter integral is often as hard as the original
problem~(\ref{eq:base}).     
Therefore it is desirable to have algorithms which are able to 
approximately compute~(\ref{eq:base}) without knowing the normalizing
constant, based solely on $n$  function values of $f$ and $\rho$. In other
terms, these functions are given by an \emph{oracle}, i.e., we assume 
that we can compute function values of $f$ and $\rho$. 

\subsubsection*{Solution operator}
\label{solop}
Assume that we are given any 
class $\fo$ of input data $(f,\rho)$ defined 
on a set $\Omega$.
We can rewrite the integral in~(\ref{eq:base}) as
\begin{equation}   \label{eq02} 
S(f, \rho) = \frac{\int f(x) \cdot \rho (x) \, \mu(dx)}{\int \rho (x)
\, \mu(dx)},\quad (f,\rho)\in\fo. 
\end{equation} 
This \emph{solution operator} is linear in $f$ but not in $\rho$. 
We discuss algorithms for the (approximate) computation of $S(f, \rho)$.
\begin{rem}
This solution operator is closely related to systems in statistical
mechanics, which obey a Boltzmann 
(or Maxwell or Gibbs) distribution, i.e., when there is a
countable number $j=1,2,\dots$ of microstates with energies, say
$E_j$,  and the overall system is distributed according to the
Boltzmann distribution, with inverse temperature $\beta$,  as
$$
P_\beta(j):= \frac{e^{-\beta E_j}}{Z_\beta},\quad j=1,2,\dots.
$$
In this case the normalizing constant $Z_\beta$ is the \emph{partition
function}, 
corresponding to $1/c$ from~(\ref{eq:base}) and $\rho^\beta(j)=
e^{-\beta E_j}$ for $j \in \N$.

In this setup, if $A$ is any global thermodynamic quantity, then its
expected value $\langle A \rangle_\beta$ is given by
$$
\langle A \rangle_\beta := \frac{1}{Z_\beta} \sum_{j} A_j e^{-\beta E_j},
$$
which can be written as $S(A,\rho^\beta)$.
Observe, however, that we use here slightly
different assumptions 
since we use the counting measure on $\N$, not a probability measure.
\end{rem}

\subsubsection*{Randomized methods}
\label{randm}

Monte Carlo methods (randomized methods) are 
important numerical tools for integration and
simulation in science and engineering, we refer to the 
recent special issue~\cite{10.1109/MCSE.2006.27}.
The Metropolis method, or more accurately, the class of 
\emph{Metropolis-Hastings algorithms} ranges among the most important 
methods in numerical analysis and scientific computation, 
see~\cite{10.1109/5992.814660,10.1109/MCSE.2006.30}.

\begin{figure}[h]
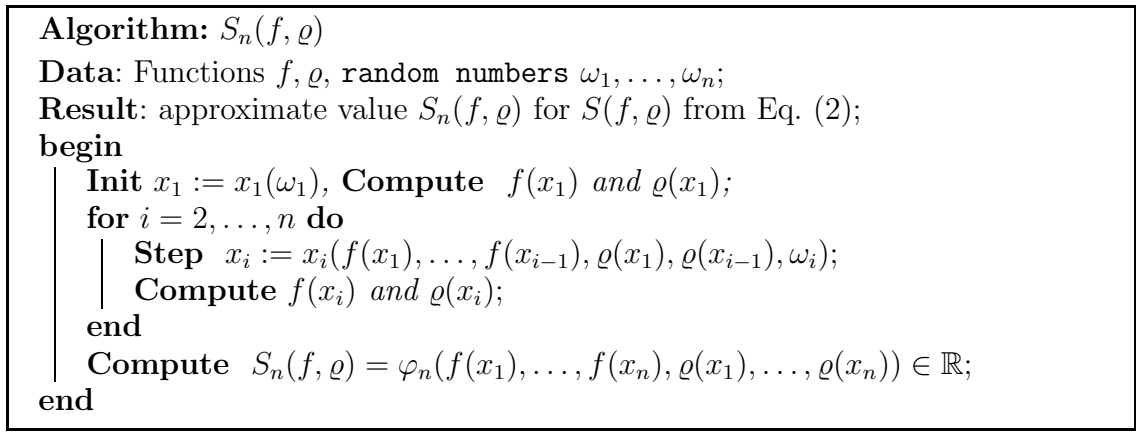

  \centering
\begin{algorithm}[H]
\SetLine
\Titleofalgo{ $S_n(f,\rho)$}
\KwData{Functions $f,\rho$, {\tt random numbers} $\omega_{1},\dots,\omega_{n}$\;}
\KwResult{approximate value $S_n(f,\rho)$ for $S(f,\rho)$ from Eq.~(\ref{eq02})\;}
\Begin{
\KwInit{$x_{1} := x_{1}(\omega_{1})$, 
\KwComp{ $f(x_1)$ and $\rho(x_1)$}\;
}

\For{$i=2,\dots,n$}
{
\KwDet{ $x_i := x_{i}(f(x_{1}),\dots,f(x_{i-1}),\rho(x_{1}),\rho(x_{i-1}),\omega_{i})$}\;
\KwComp{$f(x_i)$ and $\rho (x_i)$}\; 
}
\KwAvg{ $S_n(f,\rho)= \phi_n(f(x_{1}),\dots,f(x_{n}),\rho(x_{1}),\dots,
\rho(x_{n}))\in \R $}\;
}
\end{algorithm}  
  \caption{Generic Monte Carlo algorithm based on $n$ values of 
 $f$ and $\rho$. The final {\bf Compute} may use any mapping $\phi_n : \R^{2n} \to \R$.
%  random numbers
%  $\omega_{1},\dots,\omega_{n}$.
% Monte Carlo algorithms differ by choosing~\KwDet and~\KwAvg{} in
% different ways.
}
  \label{fig:gene}
\end{figure}

Here we consider randomized methods $S_n$ that use $n$ function 
evaluations of $f$ and $\rho$. Hence $S_n$ is of the form as exhibited
in Figure~\ref{fig:gene}. 
\SetKw{KwInit}{Init}
\SetKw{KwAvg}{Compute}
\SetKw{KwDet}{Step}
\SetKw{KwCh}{Choose}
\SetKw{KwComp}{Compute}
\restylealgo{boxed}

In all steps, random number generators may be used to determine the
consecutive node.
If the nodes $x_i$ from \KwDet
do not depend on previously computed
values of $f(x_1), \dots ,f(x_{i-1})$ and
$\rho(x_1), \dots , \rho(x_{i-1})$, then the algorithm is called 
\emph{non-adaptive}, otherwise it is called \emph{adaptive}. 
Specifically we analyze the 
procedures $\vtn$ and $\vtm$, introduced in~(\ref{eq:vtn})
and~(\ref{eq:met}) below.
\begin{rem}
The notion of \emph{adaption} which is used here differs from the
one recently used to introduce~\emph{adaptive MCMC}, see
e.g.~\cite{MR2260070,MR2172842}.
%   By a non-adaptive algorithm we mean an algorithm of the form 
%   $x_i = x_i (\omega_i)$, i.e., the node $x_i$ does \emph{not} 
%   depend on the (already computed) values 
%   $f(x_1), \dots , f(x_{i-1}), \rho(x_1), \dots , \rho (x_{i-1})$. 
%   
%   All other algorithms are called adaptive. 
The Metropolis algorithm which is used in this paper is based 
on a 
\emph{homogeneous} Markov chain, in our notation this is still
an adaptive algorithm since the used nodes $x_i$ depend on $\rho$. 
%  but the kernel of which
%  may depend on the specific target distribution, as this is the case
%  for the Metropolis sampler, see~\S~\ref{sec:metro-loc}. 
Hence we use the concept of adaptivity from numerical analysis 
and information-based complexity, see~\cite{MR1408328}. 
\end{rem}

For details on the model of computation we 
refer to~\cite{NOV,MR1319050,IBC}. 
Here we only mention the following: 
We use the real number model and assume that $f$ and $\rho$ 
are given by an oracle for function values. 
Our lower bounds hold under very general assumptions 
concerning the available random number generator.\footnote{Observe,
however, that we cannot use a random number generator 
for the ``target distribution'' 
$\mu_\rho=\rho \cdot \mu / \Vert \rho \Vert_1$, 
since $\rho$ is part of the input.} 

For the upper bounds we only study two algorithms 
in this paper, described in~(\ref{eq:vtn}) and (\ref{eq:met}),
below. Specifically we shall deal with the (non-adaptive)~\emph{simple Monte Carlo
  method} and a specific  (adaptive)~\emph{Metropolis--Hastings method}.
The former can only be applied if a random 
number generator for $\mu$ on $\Omega$ is available. 
Thus there
are natural situations when this method cannot be used.
% If applicable, then the subroutine \KwDet~ in Algorithm~\ref{fig:algorithm} chooses a random number
% according to $\mu$, independently in each step. 
The latter will be based on a suitable
ball walk. Hence we need a random number generator 
for the uniform distribution on a (Euclidean) ball.
Thus the Metropolis Hastings methods 
can also be applied when a random 
number generator for $\mu$ on $\Omega$ is not available.
Instead, we need a ``membership
oracle'' for $\Omega$: On input $x \in \R^d$ this oracle can
decide with cost 1 whether $x \in \Omega$ or not. 

% The detailed description of the \emph{adaptive} algorithm is postponed to~\S~\ref{sec:metro-loc}.
\subsubsection*{Error criterion}
\label{sec:error}
We are interested in error bounds 
uniformly for classes~$\fo$ of input data. If $S_{n}$ is any method
that uses (at most) $n$ values of $f$ and $\rho$ 
then the (individual) error for the 
problem instance~$(f, \rho)\in\fo$ is given by
\begin{equation*}
%\label{eq:mcerr}
e(S_n, (f,\rho))= \lr{\expect\abs{S(f,\rho) -
      S_n (f,\rho)}^{2}}^{1/2},
\end{equation*}
where $\expect$ means the expectation. 
The overall (or worst case)  error on the class $\fo$ is
\begin{equation*}
%\label{eq:mcerfc}
e(S_n, \fo)= \sup_{(f,\rho)\in\fo} 
e(S_n , (f,\rho)).
\end{equation*}
The complexity of the problem is given by 
the error of the best algorithm, hence we let
\begin{equation*}
%\label{eq05}
e_n (\fo) := \inf_{S_n} 
e(S_n, \fo). 
\end{equation*}
The classes~$\fo$ under consideration will always contain constant
densities~$\rho = c > 0$ and all $f$ with $\Vert f \Vert_\infty
\le 1$, hence 
$$
\mathcal F_1 (\Omega) :=\set{(f,\rho),\ \abs{f(x)}\leq 1,\
x\in\Omega, \text{ and }\ \rho = c} \subset \fo.
$$
On this class the problem~(\ref{eq02}) reduces to the classical
integration problem for uniformly bounded functions, and it is well
known that the error of any Monte Carlo method can decrease at a rate
$n^{-1/2}$, at most. Precisely, it holds true that 
$$
e_{n}(\mathcal
F_{1}(\Omega))= \frac{1}{1 + \sqrt n},
$$
if the probability~$\mu$ is non-atomic, see~\cite{olm}.
On the other hand we will only consider $(f, \rho)$ with 
$S(f, \rho) \in [-1, 1]$, hence the trivial algorithm 
$S_0=0$ always has error 1. 

For the classes $\fco$ 
and $\fad(\Omega)$,  which will be 
introduced in Section~\ref{sec:m+c},
we easily obtain the optimal order 
$e_{n}(\fo) \asymp n^{-1/2}$. 
We will analyze how $e_n(\fo)$ 
depends on the parameters 
$C$ and $\alpha$, in case $\fo:=\fco$ or
$\fo:=\fad(\Omega)$, respectively. 

We discuss some of our subsequent results and provide a short
outline. 
In Section~\ref{sec:m+c} we shall specify the methods and classes of input
data to be analyzed.
The classes $\fco$,
analyzed first in Section~\ref{s2},  contain all densities $\rho$ with 
$\sup \rho / \inf \rho \le C$. In 
typical applications we may face $C=10^{20}$. 
Then  we cannot decrease the error of optimal 
methods from 1 to $0.7$ even with sample 
size $n=10^{15}$, see Theorem 1 for more details. 
Hence the classes $\fco$ are so large that no
algorithm, deterministic or Monte Carlo, 
adaptive or non-adaptive, can provide an acceptable 
error. We also prove that the simple (non-adaptive) Monte Carlo method is almost 
optimal, no sophisticated  Markov chain Monte Carlo method can help. 

Thus we face the question whether adaptive algorithms, 
such as the Metropolis algorithm, 
help significantly on ``suitable and interesting'' subclasses of $\fco$. 
We give a positive answer for the classes 
$\fad(\Omega)$, analyzed in Section~\ref{s3}.  Here we assume that 
$\Omega \subset \R^d$ is a convex body, and that $\mu$ is the normalized Lebesgue
measure~$\muo$ on $\Omega$.  
The class~$\fad(\Omega)$ contains logconcave densities, 
where $\a$ is the Lipschitz constant 
of $\log \rho$.
We shall establish in \S~\ref{sec:non} that 
all non-adaptive methods (such as the simple Monte
Carlo method) suffer from the curse of dimension,
i.e., % for non-adaptive methods
we get similar lower bounds as for the classes $\fco$. 
 However, in \S~\ref{sec:metro-loc} we shall design and analyze
 specific (adaptive) Metropolis algorithms that are based on some
 underlying ball walks, tuned to the class parameters%  as these are the
% spacial dimension $d$ and the Lipschitz constant $\a$
. Using such algorithms we can 
break the curse of dimension by adaption. The main error estimate for
this algorithm is given in Theorem~\ref{th5}, and we conclude 
this study with further discussion in the final Section~\ref{sec:sum}.

\section{Specific methods and classes of input}
\label{sec:m+c}
We consider the approximate computation of $S(f,\rho)$
for large classes of input data. 
Since with deterministic algorithms one cannot %E  substantially 
improve 
the trivial zero algorithm (with error 1), 
we study randomized or Monte Carlo algorithms.

\subsection*{The methods}
The Monte Carlo methods under consideration  fit the schematic view from
Figure~\ref{fig:gene}.

\subsubsection*{{Simple Monte Carlo}}
\label{sec:simp}
Here the random numbers
$\omega_{1},\dots,\omega_{n}$ are identically and independently
distributed according to~$\mu$, and the routine~\KwDet chooses
$X_{i}:= \omega_{i}$. 
The final routine~\KwAvg is the quotient of the sample means of
the computed function values
\begin{equation}\label{eq:vtn}
\vtn(f,\rho):= \frac{\sum_{j=1}^n f(X_j)\rho(X_j)}{\sum_{j=1}^n\rho(X_j)}. 
\end{equation}
\subsubsection*{{Metropolis-Hastings method}}
\label{sec:mh}
This describes a class of (adaptive) Monte Carlo  methods which are based
on the ingenious idea to construct in \KwDet a Markov chain having
\begin{equation}  \label{mur} 
\mur := \frac{\rho \cdot \mu}{\int\rho(x)\, \mu(dx)}
\end{equation} 
as invariant distribution without knowing the normalization. 
Thus, if $(X_1,X_2,\dots,X_n)$ is a
trajectory of such a Markov chain, then we let \KwAvg be given as
\begin{equation}
  \label{eq:met}
  \vtm(f,\rho):= \frac{1}{n}  \sum_{j=1}^n f(X_j).
\end{equation}
Hence we use $n$ steps of the Markov chain, the number of needed 
(different) 
function values of $\rho$ and $f$ might be smaller. 
We will further specify the Metropolis-Hastings algorithm for the
problem at hand in \S~\ref{sec:metro-loc}, see Figures 2 and 3 
for a schematic presentation and Theorem 5 for the choice of $\delta$. 
%E  In diesem Bereich ein paar kleine Aenderungen. 
Both Monte Carlo methods construct Markov chains,  i.e., the point 
$x_i$ depends on $x_{i-1}$ and $\rho (x_{i-1})$, only. This trivially holds true
for simple Monte Carlo, since $x_i$ does not at all depend on 
earlier computed function values. 

\begin{rem}
Comparisons of different Monte Carlo methods for problems similar
to~(\ref{eq02}) are frequently met in the literature. We
mention~\cite{B/D06} with a comparison 
of \emph{Metropolis algorithms} and
\emph{importance sampling}, where an error expansion at any instance
$(f,\rho)$ is given in terms of certain auto-correlations. The simple
Monte Carlo method, as introduced below, is also studied there as
$\tilde\mu_{I}$ for $\rho   = 1$.
\end{rem}
The (point-wise almost sure) convergence of both 
methods $\vtn% (f,\rho)
$ and
$\vtm% (f,\rho)
$, as $n\to\infty$,  is ensured by corresponding
ergodic theorems, see~\cite{MR797411}. But, as outlined above, we are
interested in the uniform error on  relatively large~\emph{problem classes}. 
\subsection*{The classes}
Here we formally describe the classes of input under consideration.

\subsubsection*{{ The class $\fco$}}
\label{sec:classfc}

%In Section~3 we assume that 
Let $\mu$ be an arbitrary probability 
measure on a set $\Omega$ and consider the set
$$
\fco = \{ (f, \rho) \mid 
\Vert f \Vert_\infty \le 1, \ 
\rho >0, \
\frac{\rho(x)}{\rho(y)}  \le C,\ x,y\in\Omega \}. 
$$
% $$
% \fco = \{ (f, \rho) \mid 
% \Vert f \Vert_\infty \le 1, \ 
% \rho >0, \
% \frac{\sup \rho}{\inf \rho}  \le C \}. 
% $$
Note that necessarily $C\geq 1$. If $C=1$ then $\rho$ is constant and
we almost face the ordinary integration problem, since 
$\rho$ can be recovered with only one function value. 

In many applications the constant $C$ is huge and we will establish
that the complexity of the problem (the cost of an optimal
algorithm) is linear in $C$. Therefore, for large $C$, the class is
 too large. We have to look for smaller classes that 
contain many interesting pairs $(f, \rho)$ and have smaller complexity. 

\subsubsection*{The class $\mathcal \fad(\Omega)$ 
with log-concave densities}
\label{sec:classfad}

In many applications, we have a weight~$\rho$ with additional 
properties and %in Section~4 
we assume the following:
\begin{itemize} 
\item The set $\Omega\subset \R^d$ is a \emph{convex body}, that is a compact and convex set
with nonempty interior. The probability $\mu=\muo$ is the normalized Lebesgue measure 
on the set~$\Omega$. 
\item
The functions $f$ and $\rho$ are defined on $\Omega$.
\item
The weight~$\rho >0$ is log-concave, i.e., 
$$
\rho(\lambda x + (1-\lambda)y) \ge \rho(x)^\lambda \cdot 
\rho(y)^{1-\lambda}, 
$$
where $x,y \in \Omega $ and $0<\lambda <1$. 
\item
The logarithm of $\rho$ is Lipschitz, 
i.e., 
$ 
|\log\rho(x) - \log\rho(y) | \leq \alpha \Vert x-y \Vert_2
$. 
\end{itemize} 
Thus  we  consider the class of log-concave weights on
$\Omega\subset \R^{d}$ given by
\begin{equation}
\label{eq:dens-class}
\rad(\Omega)  = \{  \rho \mid 
\rho >0, \
\log\rho \text{ is concave}, \
|\log\rho(x) - \log\rho(y) | \leq \alpha \Vert x-y \Vert_2 \} . 
\end{equation}

We study the following class $\fad(\Omega)$ of problem elements, 
\begin{equation}
  \label{eq:fad}
 \fad (\Omega)  = \set{(f, \rho) \mid 
\rho \in\rad( \Omega),  \ \norm{f}{2,\rho}\le 1 } ,
\end{equation}
where $\Vert \cdot \Vert_{2,\rho}$ is the 
$L_2$-norm with respect to the probability measure $\mur$, 
see~\eqref{mur}. 
In some places we restrict our study to the (Euclidean) unit ball, i.e.,  
$\Omega:= \ball \subset \R^d$. 

\begin{rem}
Let $\radc (\Omega)$ be the class of weight functions that 
belong to $\fco$. Then 
$\rad(\Omega) \subset \radc (\Omega)$ 
if $C = e^{\alpha
D}$, where $D$ is the diameter of $\Omega$. 
Thus large $\a$ correspond to ``exponentially large'' 
values of $C$. However,
the densities from the class
$\rad(\Omega)$ have some extra (local) properties: they are log-concave 
and Lipschitz continuous. 
These properties can be used for the construction of fast 
adaptive methods, via rapidly mixing Markov chains. 
\end{rem}

\section{Analysis for $\fco$} \label{s2} 

We assume that $\Omega$ is an arbitrary set and $\mu$ 
is a probability measure on $\Omega$, 
and that the functions~$f$ and $\rho$ are defined on $\Omega$. 

In the applications, the constant $C$ might be very large, 
something like $C=10^{20}$ is a realistic assumption. 
Therefore we want to know how the complexity (the cost of 
optimal algorithms) depends on $C$. 
Observe that the problem is correctly normalized or scaled such that 
$
S(\fco) = [-1, 1] ,
$
for any $C \ge 1$. 
We will prove that the complexity of the problem % , i.e., 
% the cost of optimal algorithms,
is linear in $C$, and hence
there is no way to solve the problem if $C$ is really huge. 
% This problem class is simply too large.
We start with establishing a lower bound and then show that simple
Monte Carlo achieves this error up to a constant. 

\subsection{Lower Bounds} 

Here we prove lower bounds for all
(adaptive or non-adaptive) methods that use $n$ evaluations 
of $f$ and $\rho$. We use the technique of Bahvalov, i.e., 
we study the average error 
of deterministic algorithms with respect to certain discrete measures 
on $\fco$. 
\begin{theorem} \label{thm:lbfc}
Assume that we can partition  $\Omega$ into $2n$ disjoint sets with
equal measure (equal to $1/2n$). 
Then for any Monte Carlo method $S_n$ that uses $n$ values of 
$f$ and $\rho$ we have the lower bound
\begin{equation}
 \label{eq:2nc} 
e(S_n,\fco) \ge\frac 1 6 \sqrt 2 
\begin{cases}
\sqrt{\frac{C}{2n}}, &  2n\geq C - 1, \\
\frac{3 C}{C+2n-1}, & 2n < C -1.
\end{cases}
\end{equation}
\end{theorem} 
The lower bound will be obtained in two steps.
\begin{enumerate}
\item We first reduce the error analysis for Monte Carlo sampling to
  the average case error analysis with respect to a certain prior
  probability on the class $\fco$. 
  This approach is due to Bahvalov, see~\cite{Bachvalov}.
\item For the chosen prior the average case analysis can be carried
  out explicitly and will thus yield a lower bound. 
\end{enumerate}
To construct the prior let $m:=2n$ and  $\Omega_{1},\dots,\Omega_{m}$
the partition into sets of equal probability, and $\chi_{\Omega_{j}}$
the corresponding characteristic functions. Furthermore, let 
$$
l:=
\begin{cases}
 \lceil \frac{m}{C-1}\rceil, &  m\geq C -1,\\
1,&\text{ else.} 
\end{cases}
$$ 
Denote $J_{l}^{m}$ the set of
all subsets of $\set{1,\dots,m}$ of cardinality equal to $l$, and
$\mu_{m,l}$ the equi-distribution on $J_{l}^{m}$, while $\expect_{m,l}$ denotes the expectation with 
respect to the prior $\mu_{m,l}$. Let
$(\e_{1},\dots,\e_{m})$ be independent and identically
distributed with $P(\e_{j}=-1)= P(\e_{j}=1)=1/2,\ j=1,\dots,m$.
The overall prior is the product probability on $J_{l}^{m}\times
\set{\pm 1}^{m}$.
For any realization $\om=(I,\e_{1},\dots,\e_{m})$ we assign
$$
f_{\om}:= \sum_{j\in I} \e_{j}\chi_{\Omega_{j}}\quad \text{and}\quad
\rho_{\om}:= C \sum_{j\in I}\chi_{\Omega_{j}} + \sum_{j\not\in I}\chi_{\Omega_{j}} .
$$
The following observation is useful.
\begin{lemma}\label{lem:eml}
For any subset $N\subset\set{1,\dots,m}$ of cardinality at most $n$ it holds
$$
\expect_{m,l}\#(I\setminus N)\geq \frac l 2.
$$
\end{lemma}
\begin{proof}
  Clearly, for any fixed $k\in\set{1,\dots,m}$ we have
  $\mu_{m,l}(k\in I)=l/m$, thus
$$
\expect_{m,l}\#(I \setminus N) = \sum_{r\in N^{c}} \expect_{m,l}\chi_{I}(r) =
\#(N^{c})\frac l m\geq \frac l 2,
$$
where we denoted by $N^{c}$ the complement of $N$.

\end{proof}
\begin{proof}[Proof of Theorem~\ref{thm:lbfc}]
Given the above prior let us denote 
\begin{equation}
  \label{eq:errmfl}
  e^{avg}_{n}(\fco):= \inf_{q}\lr{\expect_{m,l}\expect_{\e}\abs{S(f,\rho) - q(f,\rho)}^{2}}^{1/2},
\end{equation}
where the $\inf$ is taken with respect to any
(possibly adaptive) deterministic algorithm
which uses at most $n$ values from $f$ and $\rho$.

For any Monte Carlo method $S_n$  we have, using Bahvalov's argument~\cite{Bachvalov}, the relation
\begin{equation}
  \label{eq:mc2avg}
 e(S_{n},\fco) \geq e^{avg}_{n}(\fco).
\end{equation}
We provide a lower bound for $e^{avg}_{n}(\fco)^{2}$. 
To this end note that for each realization $(f_{\om},\rho_{\om})$ the
integral $\int \rho_{\om} \;d\mu$ is constant.
In the first case $m\geq C -1$, 
and we can bound the integral by the choice of $l$ as
\begin{equation}
  \label{eq:intrho}
  c_{m,l}:= \int \rho_{\om}(x)\; \mu(d x)= \frac 1 m \lr{l C +
    (m-l)1} \leq 3.
\end{equation}
In the other case  $m <  C -1$, we obtain~$c_{m,1}= (C - 1 + m)/m$.
Now, to analyze the average case error, let $q_{n}$ be any
(deterministic) method, and let us assume that it uses the set $N$ of nodes. 
We have the decomposition
$$
S(f_{\om},\rho_{\om}) -
q_{n}(f_{\om},\rho_{\om})=  \lr{\frac{C}{m c_{m,l}} \sum_{j\in
    I\setminus N} \e_{j}}
 - \lr{\frac{C}{m c_{m,l}} 
 \sum_{j\in I\cap N} \e_{j} - q_{n}(f_{\om},\rho_{\om})}.
$$
Given $I$, 
the random variables in the brackets
are conditionally independent, thus uncorrelated.
Hence we conclude that
\begin{align*}
  \expect_{m,l}\expect_{\e}\abs{S(f_{\om},\rho_{\om}) -
q_{n}(f_{\om},\rho_{\om})}^{2}
& \geq \expect_{m,l}\expect_{\e}\abs{\frac{C}{m c_{m,l}} \sum_{j\in
    I\setminus N} \e_{j} }^{2}\\
& = \frac{C^{2}}{m^{2} c_{m,l}^{2}}\expect_{m,l}\#(J\setminus N)\geq
\frac{C^{2} l}{2 m^{2} c_{m,l}^{2}},
\end{align*}
by Lemma~\ref{lem:eml}.
% \begin{equation*}
%   \expect_{m,l}\abs{ S(f,\rho) - q_{n}(f,\rho)}^{2}= \frac{1}{\binom m
%     l}\sum_{I\in J_{l}^{m}}
% \expect_{m,l}\lr{\abs{ \frac{C}{m c_{m,l}} 
% \sum_{j\in I} f^{m}_{j} - q_{n}(f^{m},\rho^{m})}^{2}/I},
% \end{equation*}
% where the expectation on the right is the conditional expectation,
% i.e., when $I$ is fixed. 
% This depends on the overlap between $N\subset\set{1,\dots,m}$ the
% set of nodes which is used by $q_{n}$ and
% $I$ which may vary between $0$ and $l$. Further note that, 
%  such that we can bound ($k$ being the random
% cardinality $\#(I\setminus N)$ ) 
% \begin{align*}
%   \expect_{m,l}\abs{ S(f,\rho) - q_{n}(f,\rho)}^{2}
% &\geq  \frac{1}{\binom m l}
% \sum_{I\in J_{l}^{m}} \expect_{m,l}\lr{\abs{ \frac{C}{m c_{m,l}} 
% \sum_{j\in
%       I\setminus N} f^{m}_{j}}^{2}/I}\\
% &=  \frac{C^{2}}{m^{2} c_{m,l}^{2}} 
% \sum_{k=0}^{l} k P(\# (I\setminus N)=k)\\
% &= \frac{C^{2}}{m^{2} c_{m,l}^{2}}  \sum_{k=0}^{l} k
% \frac{\binom{n}{l-k}\binom{m-n}{k}}{\binom m l}=  
% \frac{C^{2}}{m^{2}c_{m,l}^{2}} \frac{(m - n) l}{m},
% \end{align*}
% where we used  the definition of the binomials 
% to evaluate the sum on the
% right. %  as 
% $$
%  \sum_{k=0}^{l} k
% \frac{\binom{n}{l-k}\binom{m-n}{k}}{\binom m l} = \frac{(m - n) l}{m}.
% $$
% Overall we obtain
% \begin{equation}
%   \label{eq:finbound}
%  \expect_{m,l}\abs{ S(f,\rho) - q_{n}(f,\rho)}^{2} \geq
%  \frac{C^{2}}{m^{2}c_{m,l}^{2}} \frac{(m - n) l}{m}= \frac{C^{2} l}{2
%    c_{m,l}^{2} m^{2}}.  
% \end{equation}
In the case $m\geq C -1 $ we obtain $l\geq m/C$ and  have
$c_{m,l}\leq 3$, such that %we finally obtain
$$
\expect_{m,l}\abs{ S(f,\rho) - q_{n}(f,\rho)}^{2} \geq \frac{C}{36 n},
$$
which in turn yields the first case bound in~(\ref{eq:2nc}).
In the other case~$m <  C -1$ the value of $l=1$ yields the second
bound in~(\ref{eq:2nc}).
\end{proof} 
\subsection{The error of the simple Monte Carlo method} 
\label{sec:simple}

The direct approach to evaluate~(\ref{eq:base}) would be to use the
method~$\vtn$ from~(\ref{eq:vtn}).
We will prove an upper bound for the error of this method, and 
we start with the following 
\begin{lemma}\label{lem:rho}
  If the function $\rho$ obeys the requirements in~$\fco$, then 
  \begin{enumerate}
  \item $0< \inf_{x\in\Omega}\rho(x)\leq
    \sup_{x\in\Omega}\rho(x)<\infty$.
\item For every probability measure $\mu$ on $\Omega$ we have
$\norm{\rho}{2,\mu}\leq \sqrt C\norm{\rho}{1,\mu} $.
  \end{enumerate}
\end{lemma}
\begin{proof}
  To prove the first assertion, fix any $y_{0}\in\Omega$. Then the
  assumption on $\rho$ yields $\rho(x)\leq C \rho(y_{0})$, and
  reversing the roles of $x$ and $y$ also the lower bound.
Now both, the assumption on $\rho$ as well as the 
second assertion,  are invariant with respect to multiplication
of $\rho$ by a constant. In the light of the first assertion we may
and do assume that $1\leq\rho(x)\leq C,\ x\in\Omega$, and  we derive,
using $ 1 \leq \int_{\Omega}\rho(x)\; \mu(dx)$, that
$$
\int_{\Omega}\rho^{2}(x)\; \mu(dx)\leq C \int_{\Omega}\rho(x)\;
\mu(dx) \leq C \lr{\int_{\Omega}\rho(x)\; \mu(dx)}^{2},
$$
completing the proof of the second assertion and of the lemma.
\end{proof}
We turn to the bound for the simple Monte Carlo method.
\begin{theorem}
For all $n\in\N$ we have
\begin{equation}
    \label{eq:thm1}
e(\vtn,\fco)\leq 2\, \min\set{1,  \sqrt{\frac{2C}{n}}} .  
  \end{equation}
\end{theorem}
\begin{proof}
The upper bound~$2$ is trivial, it even holds deterministically. 
  Fix any pair $(f,\rho)$ of input. For any sample
  $\lr{X_1,\dots,X_n}$ and function $g$ we denote the sample 
mean by $\vt(g):= 1/n\sum_{j=1}^n g(X_j)$. 
It is well known that $e(\vt,g)\leq \norm{g}{2}/\sqrt n$. 
With this notation we can
bound
  \begin{align*}
    &\abs{S(f,\rho) - \vtn(f,\rho)}\leq \abs{S(f,\rho) -
      \frac{\vt(f\rho)}{\int \rho(x)\mu(dx)}}+  
\abs{\frac{\vt(f\rho)}{\int \rho(x)\mu(dx)} -
  \frac{\vt(f\rho)}{\vt(\rho)}}\\ 
&\leq \frac{1}{\norm{\rho}{1}}\lr{\abs{\int
  f(x)\rho(x)\mu(dx)-\vt(f \rho) } 
+ \abs{\frac{\vt(f\rho)}{\vt(\rho)}}
\abs{\int \rho(x)\mu(dx) - \vt(\rho)}}\\
&\leq  \frac{1}{\norm{\rho}{1}}\lr{\abs{\int
  f(x)\rho(x)\mu(dx)-\vt(f \rho) } 
+\norm{f}{\infty}
\abs{\int \rho(x)\mu(dx) - \vt(\rho)}},
  \end{align*}
where we used 
$
\abs{\vt(f\rho)/{\vt(\rho)}}\leq \norm{f}{\infty},
$
which holds true since the enumerator and 
denominator use the same sample.
This yields the following error bound
\begin{align*}
  e(\vtn,(f,\rho))&\leq   \frac{\sqrt 2}{\norm{\rho}{1}}
\lr{ e(\vt,f\rho) + \norm{f}{\infty}e(\vt,\rho)}\\
&\leq \frac{\sqrt 2}{\norm{\rho}{1}\sqrt n}\lr{\norm{f\rho}{2} +
  \norm{f}{\infty}\norm{\rho}{2}}\leq \frac{2\sqrt 2 \norm{f}{\infty}}{\sqrt
  n}\frac{\norm{\rho}{2}}{\norm{\rho}{1}}
  \leq \frac{2\sqrt{2C}}{\sqrt  n},
\end{align*}
where we use Lemma~\ref{lem:rho}. Taking  the supremum over $(f,\rho)\in\fco$
allows to complete the proof. 
\end{proof}

\section{Analysis for $\fad(\Omega)$} \label{s3} 

In this section we impose restrictions on  the input data, in
particular on the density,  in order to improve the complexity. This
class is still large enough to contain many important situations.
Monte Carlo methods for problems when the target (invariant)
distribution is log-concave proved to be important in many studies, we
refer to~\cite{MR1284987}. One of the main intrinsic features of such
classes of distributions are \emph{isoperimetric inequalities},
see~\cite{103439,MR1318794}, which will also be used here in the form
as used in~\cite{MR2178341}.
Recall that here we always require that $\Omega\subset \R^{d}$ is a
convex body, as introduced in Section~\ref{sec:classfad}.

% We always assume the following. 
% The functions $f$ and $\rho$ are defined on 
% a compact and convex set $\Omega \subset \R^d$ 
% with nonempty interior and  
% $\mu=\muo$ is the normalized Lebesgue measure 
% on the set~$\Omega$. 
% We  consider the class of log-concave weights on
% $\Omega$  given by
% $$ 
% \rad(\Omega)  = \{  \rho \mid 
% \rho >0, \
% \log\rho \text{ is concave}, \
% |\log\rho(x) - \log\rho(y) | \leq \alpha \Vert x-y \Vert_2 \} . 
% $$

% We study the class $\fad(\Omega)$, given by 
% $$
% \fad (\Omega)  = \set{(f, \rho) \mid 
% \rho \in\rad( \Omega),  \ \norm{f}{2,\rho}\le 1 } ,
% $$
% where $\Vert \cdot \Vert_{2,\rho}$ is the 
% $L_2$-norm with respect to the probability measure $\mur$, 
% see~\eqref{mur}. 
% In particular, we study the (Euclidean) unit ball, i.e.,  
% $\Omega:= \ball \subset \R^d$. 

We start with a lower bound for all non-adaptive algorithms to exhibit
that simple Monte Carlo cannot take into account the additional
structure of the underlying class of input data and adaptive methods
should be used. This bound, together with Theorem~\ref{th5}, will show 
that adaptive methods can outperform any 
non-adaptive method, if we consider $S$ on $\fad (\ball)$. 
Indeed, we also show that specific Metropolis
algorithms, based on local underlying Markov chains are suited for
this problem class.

\subsection{A lower bound for non-adaptive methods}
\label{sec:non}

Here we prove a lower bound for all non-adaptive methods 
(hence in particular for the simple Monte Carlo method) 
for the problem on the classes~$\fad(\Omega)$. 
Again, this lower bound will use Bahvalov's technique.

We start with a result on sphere packings. 
The Minkowski-Hlawka theorem,  see~\cite{MR0172183}, 
says that the density of the densest sphere packing in $\R^d$ 
ist at least $\zeta (d) \cdot 2^{1-d}\ge 2^{1-d}$. 
It is also known, see \cite{Hlawka}, that the density 
(by definition of the whole $\R^d$) can be replaced by the density within 
a convex body $\Omega$, as long as the radius $r$ of the 
spheres tends to zero. Hence we obtain the following result. 

\begin{lemma}
\label{lem:MHT}
There is $n_{\Omega}\in\N$ such that for all $m\geq n_{\Omega}$ there are points
$y_{1},\dots,y_{m}\in\Omega$ such that with 
$$
r:=r(\Omega,m):= 2^{-1} m^{-1/d} \left( \frac{\vol (\Omega)}
{\vol (\ball)}\right)^{1/d}
$$
the closed balls $B_{i}:= B(y_{i},r)\subset \Omega$ 
are disjoint.
\end{lemma}

Our construction will use such points $y_{1},\dots,y_{m}\in\Omega$ and
the corresponding balls $B_{1},\dots,B_{m}$ as follows.

For $i\in\set{1,\dots,m}$ we assign 
\begin{align*}
 \rho_{i}(y)&:= c_i \exp\lr{-\alpha\norm{y - y_{i}}{2}},\quad
 y\in\Omega  \quad\text{and}\\ 
f_{i}(y)&:= \tilde c_i  \chi_{B_{i}}(y),\quad y\in\Omega ,
\end{align*}
with constants $c_i$ and $\tilde c_i$ chosen such that 
\begin{alignat*}{2}
1&= \int_{\Omega } \rho_i(y) \, dy &= 
c_i \int_{\Omega } \exp(- \a \norm{y - y_i}{}) dy\quad \text{and}\\
1&=\norm{f_i}{2,\rho_i} &= \tilde c_i^2 c_i \int_{B_i} \exp(- \a
\norm{y - y_i}{})\, dy.  
\end{alignat*}
The corresponding values of the mapping $S$ are computed as
\begin{align}\label{eq:slb}
  \begin{split}
S(f_i,\rho_i) &= \int_{\Omega } f_i \rho_i\, dy = \tilde c_i c_i
\int_{B_i} \exp(- \a \norm{y - y_i}{})\, dy\\
& = \lr{ c_i \int_{B_i} \exp(- \a \norm{y - y_i}{}) dy}^{1/2}=
 \lr{ c_i \int_{B(0,r)} \exp(- \a \norm{y}{}) \, dy}^{1/2}\\
&= \lr{\frac{\int_{B(0,r)} \exp(- \a \norm{y}{}) \, dy}{\int_{\Omega} 
\exp(- \a \norm{y - y_i}{})\,dy}}^{1/2}.
  \end{split}
\end{align}
Again we turn to the average case setting, this time with
 probability measure $\mu^{2n}$ being the equidistribution on the set 
$$
\mathcal F^{2n}:= \set{ \lr{\e_i f_i,\rho_i},\quad i=1,\dots,2n,\
  \e_i=\pm 1}\subset \fad(\Omega ).
$$
Similar to~(\ref{eq:mc2avg}) we have for any non-adaptive Monte Carlo
method $S_n(f,\rho)$ the relation 
$$
e(S_n,\fad(\Omega ))\geq
\min\set{ e^{avg}(q_n,\mu^{2n}),\quad q_n \text{ is 
deterministic and non-adaptive}},
$$
where $e^{avg}(q_n,\mu^{2n})$ denotes the average case error of the
deterministic non-adaptive method $q_n$ with respect to the
probability $\mu^{2n}$.
Thus let~ $q_n$ be any non-adaptive 
(deterministic) algorithm for $S$ on the 
class $\fad (\Omega )$ that uses at most $n$ values.

The average case error can then be bounded from below as
\begin{align*}
\expect_{\mu^{2n}}\abs{S(f,\rho) - q_n(f,\rho)}^2&=
\frac{1}{2n}\sum_{i=1}^{2n} 
\expect_{\e}\abs{S(\e_i f_i,\rho_i) - q_n(\e_i
  f_i,\rho_i) }^2\\
&\geq \frac 1 2 \min_{i=1,\dots,2n}\expect_{\e}\abs{S(\e_i f_i,\rho_i)
}^2 \geq  \frac 1 2 \min_{i=1,\dots,2n}S(f_i,\rho_i)^2.
\end{align*}
Above, $\expect_{\e}$ denotes the expecation with respect to the
independent random variables $\e_{i}=\pm 1$.
Together with~(\ref{eq:slb}) we obtain
$$
e(S_n,\fad(\Omega))\geq \frac 1 2 \sqrt 2\,
\min_{i=1,\dots,2n}\lr{\frac{\int_{B(0,r)} \exp(- \a \norm{y}{}) \,
    dy}{\int_{\Omega} \exp(- \a \norm{y - y_i}{})\,dy}}^{1/2}. 
$$
We bound the enumerator from below and the denominator from
above.
For $\alpha r\leq \log 2$ we can bound 
$$
\int_{B(0,r)} \exp(- \a \norm{y}{}) \, dy\geq \frac 1 2
  \vol(B(0,r))= \frac 1 2 r^d \vol(\ball).
$$
For the denominator we have  %  \fix{alpha gross klein} 
%   letting temporarily $\bar\a:=
%   \max\set{\a,1}$, that 
\begin{align*}
\int_{\Omega} \exp(- \a \norm{y - y_i}{})\,dy &\leq \int_{\R^d} \exp(-
\a \norm{y - y_i}{})\,dy \\
& ={\a}^{-d} \int_{\R^d} \exp(-\norm{y}{})\,dy=
{\a}^{-d} \Gamma(d)\vol{\partial \ball},
\end{align*}
such that we finally obtain, using the well known formula
$\vol(\partial \ball) = d \vol(\ball)$, that
$$
e(S_n,\fad(\Omega))\geq  \frac 1 2 \sqrt 2\, \lr{\frac{{\a}^d
    r^d}{2 d!}}^{1/2} = \frac 1 2 \lr{\frac{{\a}^d
    r^d}{d!}}^{1/2}.
$$
Using the value for $r=r(\Omega ,2n)$ from Lemma~\ref{lem:MHT} we end up
with
\begin{theorem} 
Assume that $S_n$ is any non-adaptive Monte Carlo method for 
the class $\fad (\Omega )$. Then, with~$ n_\Omega $ from Lemma~\ref{lem:MHT},
we have for all 
$$ 
2n \ge \max\set{n_\Omega ,\lr{\a/{\log 4}}^d \cdot 
\frac{\vol \Omega}{\vol \ball}}
$$ 
that
\begin{equation}  \label{lo9} 
e(S_n, \fad(\Omega )) \ge 
2^{-d/2-3/2} \cdot 
\left( \frac{\vol \Omega}{\vol \ball} \right)^{1/2} \cdot 
\frac{\alpha^{d/2}}{\sqrt{d!}} \ n^{-1/2} .
\end{equation} 
\end{theorem} 

\begin{rem}
For fixed $d$ this is a lower bound of the form 
$e(S_n) \ge c_\Omega \, \a^{d/2} \, n^{-1/2}$. It is interesting only 
if $\alpha$ is ``large'', otherwise the already mentioned lower bound 
$(1+ \sqrt{n})^{-1}$ is better. 

We stress that in the above reasoning we essentially used the
non-adaptivity of the method $S_n$. Indeed, if $S_n$ were adaptive,
then by just one appropriate function 
value $\rho(x)$,  we could identify the
index $i$, since the functions $\rho_i$ are
global. Then, knowing $i$,  we could ask for the value of $\e_i$ and
would obtain the exact solution to $S(f,\rho)$ for this small class
$\mathcal F^{2n}$ for all $n \ge 2$. 
\end{rem}

\subsection{Metropolis method with local underlying walk}
\label{sec:metro-loc}

The Metropolis algorithm we consider here has a specific
routine~\KwDet in Figure~\ref{fig:gene}, whereas the
final step~\KwAvg is exactly as given in~(\ref{eq:met}). It is based on a
specific ball walk and this version is
sometimes called \emph{ball walk with
Metropolis filter}, see~\cite{MR2178341}.
Two concepts from the theory of Markov chains turn out to be
important, reversibility and uniform ergodicity. We recall these
notions briefly, see~\cite{MR1399158} for further details.
A Markov chain  $(K,\pi)$ is \emph{reversible with respect to $\pi$}, 
if for all measurable subsets $A,B\subset\Omega$ the balance
\begin{equation}\label{eq-rev}
\int_{A}K(x,B)\pi(dx)=\int_{B}K(x,A)\pi(dx)
\end{equation}
holds true. Notice that in this case necessarily $\pi$ is an invariant
distribution.
 
A Markov chain is \emph{uniformly ergodic} if there are $n_{0}\in\N$, a
constant $c>0$ and a probability measure $\nu$ on $\Omega$ such that
  \begin{equation}
    \label{eq:ueball}
K^{n_{0}}(x,A) \geq c \nu(A),
\quad \text{ for all } A\subset \Omega\text{ and } x\in\Omega.
  \end{equation}
Markov chains which are  uniformly ergodic have a unique invariant
probability distribution.

Our analysis will be based on conductance arguments and we
recall the basic notions, see~\cite{MR1025467,MR1238906}.
If $(K,\pi)$ is a Markov chain with transition kernel $K$ and
invariant distribution $\pi$ then we assign the 
\begin{enumerate}
\item 
\emph{local conductance} at $x\in\Omega$ by $l_K(x):=
  K(x,\Omega\setminus\set{x})$,
\item and the \emph{conductance} as
\begin{equation}
  \label{eq:conductance}
  \phi(K,\pi):= \inf_{0<\pi(A) <  1}\frac{\int_A K(x,A^c)
    \pi(dx)}{\min\set{\pi(A),\pi(A^c)}}, 
\end{equation}
where $A^c= \Omega \setminus A$. 
\end{enumerate}
Below we call $l>0$ a \emph{lower bound for the local conductance}, if
$l_{K}(x)\geq l$ for all $x\in\Omega$.

\subsubsection*{The ball walk and some of its properties}
\label{sec:ball}

Here we gather some properties of the ball walk, 
see~\cite{MR1238906,MR2178341},  which will serve as
ingredients for the analysis of Metropolis chains using this as the
underlying proposal. 
In particular we prove that on convex bodies in $\R^{d}$ the ball walk is 
uniformly ergodic and we bound its conductance from below, in terms
of bounds $l>0$ for the local conductance.

We abbreviate $B(0,\delta) = \delta \ball$. 
Let $Q_\delta$ be the transition 
kernel of a local random walk
having transitions within $\delta$-balls of its current position,
i.e., we let 
\begin{equation}
\label{eq:pxx}
Q_{\delta}(x,\set{x}):= 1 - \frac{\vol(B(x,\delta) 
\cap \Omega)}{\vol(\delta \ball )},
\end{equation}
and 
\begin{equation}
\label{eq:qloc}
Q_\delta(x,A):= 
\begin{cases}
\displaystyle{\frac{\vol(B(x,\delta) \cap A)}{\vol(\delta \ball )}}, 
&   A
\subset \Omega \text{ and }x \notin A, \\
Q_\delta(x,A\setminus\set{x}) +   Q_{\delta}(x,\set{x}), & 
A \subset \Omega \text{ and } x\in A.
\end{cases}
\end{equation}
Schematically, the transition kernel may be viewed 
as in Figure~\ref{fig:bbb}.

\SetKw{KwProp}{Propose:}
\SetKw{KwAcc}{Accept:}
\SetKwInOut{Input}{Input}
\SetKwInOut{Output}{Output}
\restylealgo{ruled}
\begin{figure}[h]
  \centering
\begin{procedure}[H]
\Input{current position $x$; $\delta>0$\;}
\Output{next position\;}
\KwProp{Choose $y\in B(x,\delta)$ uniformly}\;
\KwAcc{}
\eIf{$y\in\Omega$}{\Return{$y$}\;}{\Return{$x$}\;} 
  \caption{Ball-walk-step($x,\delta$)}
\end{procedure}  
  \caption{Schematic view of ball walk step}
  \label{fig:bbb}
\end{figure}
Clearly we may restrict to $\delta\leq D$, the diameter of $\Omega$.
The following observation is important and explains why we restrict
ourselves to convex bodies..
\begin{lemma}
 If $\Omega\subset \R^{d}$ is a convex body, then the ball walk
 $Q_{\delta}$ has a (non-trivial) lower bound $l>0$ for the local conductance.
\end{lemma}
\begin{proof} 
It is well-known that convex bodies satisfy the cone condition 
(see % Lemma 3 of Section 3.2 in
\cite[\S~3.2, Lemma~3]{Burenkov}).  
Therefore we obtain that for each $\delta>0$ there is $l>0$ such that
for each $x \in \Omega$ we have $l_{Q_\delta} (x) \ge l$.
% $$
% \exists \ {l > 0} \quad
% \exists \ {\delta_0>0} \quad 
% \forall \ {0<\delta < \delta_0} \quad
% \forall \ {x \in \Omega} \quad
% l_{Q_\delta} (x) \ge l .
% $$
\end{proof}
\begin{rem}
Observe however, that $l$ might be very small. 
For $\Omega=[0,1]^d$, for example, we get $l = 2^{-d}$, 
even if $\delta$ is very small. In contrast, we will 
see that a large $l$ is possible for $\Omega=B^d$ 
and $\delta \le 1/\sqrt{d+1}$, see Lemma~\ref{lem:l-bound}.   
\end{rem}
Notice that $l_{Q_{\delta}}(x)= {\vol(B(x,\delta) \cap
\Omega)}/{\vol(\delta \ball )}$, hence in the following we use the inequality
\begin{equation}
\label{eq:l-bound}
\vol(B(x,\delta)\cap \Omega)\geq l \vol(\delta\ball),  
\end{equation}
where $l>0$ is a lower bound for the local conductance 
of the ball walk. 

The following result  is~\emph{folklore}, but for a 
lack of reference we sketch a proof.

\begin{proposition}\label{prop:ueqd}
%  Let $\Omega\subset \R^{d}$ be compact.
The ball walk $Q_{\delta}$ is reversible with respect to the uniform
distribution $\muo$ and 
%  If there is a lower bound $l>0$  for the local conductance
%  of $Q_{\delta/2}$ then the ball walk $Q_{\delta}$ is 
uniformly ergodic.
%   For each $0<\delta\leq 1/\sqrt{d+1}$  the ball 
%   walk $Q_{\delta}$  is uniformly
%   ergodic and reversible. In particular there are $n_{0}\in\N$, a
%   constant $c>0$ and a probability measure $\nu$ on $\ball$ such that
%   \begin{equation}
% %    \label{eq:ueball}
% Q_{\delta}^{n_{0}}(x,A) \geq c \nu(A),
% \quad \text{ for all } A\subset \ball\text{ and } x\in\ball.
%   \end{equation}
\end{proposition}

 The crucial tool for proving this is provided by the
notion of small and petite sets, where we refer to~\cite[Sect.~5.2 \&
5.5]{Meyn-book} for details and properties. 
To this end we introduce a \emph{sampled} chain, say
$(Q_{\delta})_{a}$, where $a$ is some probability
$a=\lr{a_{0},a_{1},\dots}$ on $\set{0,1,2,\dots}$
and $(Q_{\delta})_{a}$ is defined by $(Q_{\delta})_{a}(x,C):=
\sum_{j=0}^{\infty}a_{j}Q_{\delta}^{j}(x,C)$.
%A set $C\subset \Omega$ is \emph{petite}, 
We recall that a
(measurable) subset $C\subset \Omega$ is \emph{petite} (for
$Q_{\delta}$), if there are a probability~$a$
 and a probability measure $\nu$ on
$\Omega$ such that 
\begin{equation}
\label{eq:small}
(Q_{\delta})_{a}(y,A)\geq \varepsilon \nu(A),
\quad A\subset \Omega,\ y \in C.
\end{equation}
A set $C\subset \Omega$ is \emph{small}, if the same property holds
true for some Dirac probability $a:= \delta_{n}$, such that obviously
small sets are petite.
We first show that certain balls are small.

\begin{lemma}\label{lem:small}
% Let $\delta>0$ and let $l >0$ be a lower bound for 
% the local conductance
% of the ball walk $Q_{\delta/2}$. 
The sets $ B(x,\delta/2)\cap
% If there is a is  
% a lower bound $l>0$ for the local conductance
% of the ball walk $Q_{\delta/2}$ then the sets $ B(x,\delta/2)\cap
\Omega,\ x\in\Omega$ are small for $Q_\delta$.
%  Let $\delta\leq 1/\sqrt{d+1}$ and $x\in\ball$. If $y\in
%   B(x,\delta/2)\cap \ball$ then
%   \begin{equation}
%     \label{eq:smlemma}
%     Q_{\delta}(y,A) \geq 0.3 \cdot 2^{-d} \frac{\vol(A \cap
%       B(x,\delta/2)\cap \ball)}{\vol( B(x,\delta/2)\cap \ball)},\quad
%     A\subset \ball.
% \end{equation}
%Consequently,  each set $B(x,\delta/2)\cap \ball$  is small.
\end{lemma}

\begin{proof}
First, we note that $y\in B(x,\delta/2)$ implies $B(x,\delta/2) \subset
B(y,\delta)$. Let $l>0$ be a lower bound for the local conductance of
$Q_{\delta/2}$. Using~(\ref{eq:l-bound}) 
for $Q_{\delta/2}$, we obtain for any set $A\subset \Omega$ that
\begin{align*}
  Q_{\delta}(y,A) &\geq  Q_{\delta}(y,A\setminus\set{y}) =
  \frac{\vol(B(y,\delta)\cap A)}{\vol(B(y,\delta))} \geq 2^{-d}
  \frac{\vol(B(x,\delta/2)\cap A)}{\vol(\delta/2\ball)}\\
&\geq l \cdot 2^{-d} \frac{\vol(A \cap
      B(x,\delta/2)\cap \Omega)}{\vol( B(x,\delta/2)\cap \Omega)}.
\end{align*}
Hence estimate~(\ref{eq:small}) holds true with $n_{0}:=1,\
\varepsilon:= l\cdot 2^{-d}$ and 
$$
\nu(A) := \frac{\vol(A \cap
      B(x,\delta/2)\cap \Omega)}{\vol( B(x,\delta/2)\cap \Omega)},\quad
    A\subset \Omega.
$$
This completes the proof.
\end{proof}
\begin{proof}[Proof of Proposition~\ref{prop:ueqd}]
We first prove reversibility with respect to $\muo$. 
Notice that it is enough to verify~(\ref{eq-rev}) 
for disjoint sets $A,B\subset \Omega$.
Furthermore we observe that for any pair $A,B\subset \Omega$ 
of measurable subsets the characteristic function of the set 
$$
\set{(x,y)\in\Omega\times \Omega,\quad x\in A,\ y\in B,\ \norm{x -
    y}{}\leq \delta}
$$
can equivalently be rewritten as
$$
\chi_{B}(y) \chi_{B(y,\delta)\cap A}(x)
\quad \text{or} \quad \chi_{A}(x) \chi_{B(x,\delta)\cap B}(y).
$$
Hence, letting temporarily
$c:={\vol(\Omega)\vol(\delta\ball)}$ we obtain  
\begin{align*}
  \int_{A}Q_{\delta}(x,B)\;\muo(dx)&=
 \frac 1 c \int_{A}\vol(B(x,\delta)\cap
  B)\; dx\\
&=  \frac 1 c \int_{\Omega}\int_{\Omega}\chi_{A}(x) 
\chi_{B(x,\delta)\cap B}(y)\;
dy\;dx\\
&=  \frac 1 c \int_{\Omega}\int_{\Omega}\chi_{B}(y) 
\chi_{B(y,\delta)\cap A}(x)\;
dx\;dy= \int_{B}Q_{\delta}(y,A)\;\muo(dy),
\end{align*}
proving reversibility.

By Lemma~\ref{lem:small} each set $B(x,\delta/2) \cap \Omega$ is small,
thus also petite. Petiteness is in\-heri\-ted by taking finite
unions. Since $\Omega$, being compact, can be covered by finitely many
sets  $B(x,\delta/2)\cap \Omega$, this implies that $\Omega$ is
petite. By~\cite[Thm.~16.2.2]{Meyn-book} this yields uniform
ergodicity of the ball walk % , and hence that $\Omega$ is small
(see~\cite[Thm.~16.0.2(v)]{Meyn-book}).
\end{proof}
We mention the following conductance bound  of the ball
walk, which is  a slight improvement
of~\cite[Thm.~5.2]{MR2178341}. This will be  a special case of
Theorem~\ref{thm:met-cond}, below, and we omit the proof.

\begin{proposition}\label{pro:phi}
Let $(Q_{\delta},\muo)$ be the ball walk from above,
and let $\phi(Q_{\delta},\muo)$ be its conductance. 
Let~$D$ be the diameter of $\Omega$  and 
let $l$ be a lower bound for the local conductance. Then
\begin{equation}
\label{eq:ballconductancelb}
\phi(Q_{\delta},\muo) \geq  
\sqrt{\frac \pi 2}\frac{l^{2}\delta}{8 D \sqrt{d +1}}.
%% \frac{l^{2}\delta}{16 D \sqrt d}.
\end{equation}  
\end{proposition}

The local conductance may be arbitrarily small if the domain $\Omega$
has sharp corners. 
For specific sets $\Omega$ we can explicitly provide lower bounds for
the local conductance, and this will be used in the later convergence
analysis.
In the following we mainly discuss the case $\Omega = \ball$. 

We start with a  technical result, related to the Gamma function on
$\R^+$. We use the well-known formula
\begin{equation}
  \label{eq:3}
\vol(\ball)= \pi^{d/2}/\Gamma(d/2 +1). 
\end{equation}
\begin{lemma}\label{lem:bou}
For any $z>0$ we have
\begin{equation}
  \label{eq:gamma}
  \frac{\Gamma(z+1/2)}{\Gamma(z)}\leq \sqrt z.
\end{equation}
Consequently,
\begin{equation}
  \label{eq:vol-bound}
\frac{ \vol(B^{d-1})}{\vol(\ball)}
%% \frac{\Gamma(d/2 +1 )}{\delta\sqrt\pi\Gamma((d+1)/2) }
\leq \sqrt{\frac{d+1}{2\pi}}.
\end{equation}
\end{lemma}
\begin{proof}
  By~\cite[Chapt.~VII, Eq.~(11)]{MR2013000} we know that the function
  $z\mapsto \log\Gamma(z)$ is convex for $z>0$. Thus we conclude
  \begin{align*}
    \log\Gamma(z + 1/2) 
&\leq \frac 1 2 \lr{\log\Gamma(z+1) + \log\Gamma(z)}\\
& = \frac 1 2 \lr{\log z  + 2 \log\Gamma(z)} 
= \log\sqrt z + \log\Gamma(z),
  \end{align*}
from which the proof of assertion~(\ref{eq:gamma}) can be completed.
Using the representation for the volume from~(\ref{eq:3}) and applying
the above  bound with $z:= (d+1)/2$ we obtain
$$
\frac{ \vol(B^{d-1})}{\vol(\ball)}\leq
\frac{\Gamma(d/2 +1 )}{\sqrt\pi\Gamma((d+1)/2) }
\leq \sqrt{\frac{d+1}{2\pi}},
$$
and the proof is complete.
\end{proof}
%P where we can explicitly provide lower bounds 
%  for the local conductance. 
Using Lemma~\ref{lem:bou},  we can prove the 
following lower bound for the local
conductance of the ball walk on $\ball$.

\begin{lemma}  \label{lem:l-bound}
Let $(Q_\delta,\muo)$ be the local ball walk on $\ball\subset \R^d$.
If $\delta\leq 1/\sqrt{d +1}$, then its 
local conductance obeys $l\geq 0.3$.  
\end{lemma}

\begin{proof}
The proof is based on some geometric reasoning. It is clear that the
local conductance~$l(x)$ is minimal for points $x$ at the
boundary of $\ball$, and in this case 
its value equals the portion, say $\widetilde V$,  
of the volume of $B(x,\delta)$ inside $\ball$. If $H$ is the
hyperplane at $x$ to $\ball$, then this cuts off $B(x,\delta)$
exactly one half of its volume. 
Thus we let  $Z(h)$ be the cylinder with
base being the $(d-1)$-ball around $x$ in the hyperplane $H$ of
radius $\delta$. 
Its height~$h$ is the distance of $H$ to the hyperplane
  determined by the intersection of $\ball\cap B(x,\delta)$. This
  height $h$ is exactly determined from the quotient $h/\delta =
  \delta/2$, by similarity, hence $h:= \delta^2/2$. By
  construction we have $\widetilde V \geq 1/2 -
\vol(Z(h))/\vol(B(x,\delta))$ and we can 
lower bound the local conductance $l(x)$ by
$$
l(x)\geq \frac 1 2  - \frac{\vol(Z(h))}{\vol(B(x,\delta))}.
$$
We can evaluate~$\vol(Z(h))$ as
$
\vol(Z(h)) = h \delta^{d-1} \vol(B^{d-1}),
$
and we obtain
$$
l(x)\geq \frac 1 2 - \frac{\delta^{d+1} \vol(B^{d-1})}{2 \delta^d
\vol(\ball)}= \frac 1 2 \lr{ 1 - 
\frac{\delta \vol(B^{d-1})}{\vol(\ball)}}.
$$
% We use~(\ref{eq:3}) 
% \begin{displaymath}
%   l(x)\geq \frac 1 2 \lr{1 - \frac{\delta \Gamma(\frac d 2 + 1)}{
%   \sqrt\pi  \Gamma(\frac d 2 +  \frac 1 2)}}.
% \end{displaymath}
The bound~(\ref{eq:vol-bound}) from Lemma~\ref{lem:bou} implies
$$
l(x) \geq % \frac 1 2 \lr{1 - \frac{\delta \sqrt{{(d+1)}/{2}}}{
%   \sqrt\pi}} =
\frac 1 2 \lr{1 - \frac{\delta \sqrt{{d+1}}}{  \sqrt{2
    \pi}}}.
$$
For $\delta\leq 1/(\sqrt{d+1})$ we get
$l(x) \geq 1/2( 1 - 1/\sqrt{2\pi})\geq 0.3$, completing the proof.
\end{proof}

We close this subsection with the following technical lemma, 
which  can be extracted from the unpublished
seminar note~\cite{vempala-lesson}. For the convenience of the
reader we present its proof. 
In addition we will slightly improve the statement.
\begin{lemma}%% [{\cite{vempala-lesson}}]
  \label{lem:vempala}
Let $l >  0$ be a lower bound for the local 
conductance of the ball walk $(Q_\delta,\muo)$.
For any $0<t< l$ and any set
  $A\subset \Omega$ with related sets 
  \begin{align}
A_1&:= \set{x\in A, \quad Q_\delta (x, A^c)< \frac{l -t}{2}}\subset
    A\\
  A_2 &:= \set{y\in A^c,\quad  Q_\delta(y, A)< \frac{l -t}{2}}\subset
  A^c,
  \end{align}
we have $d(A_1,A_2)>t\delta \sqrt{2 \pi/\lr{d+1}}$.
\end{lemma}
For its proof we need the following 
\begin{lemma}
Let $\delta>0$.
  If $x,y\in \R^d$ are two points with distance~$t\delta \sqrt{2
    \pi/\lr{d+1}}$ at most, then 
  \begin{equation}
    \label{eq:1}
    \vol(B(x,\delta)\cap B(y,\delta)) \geq (1 - t) \vol(\delta\ball).
  \end{equation}
\end{lemma}
\begin{proof}
Let $u:= \norm{x - y}{2}$. If $u<\delta$ then 
 the volume of the intersection of $B(x,\delta)$ and $B(y,\delta)$  is
exactly the same as the volume of the 
ball $\delta\ball$ minus the volume of the
 middle slice with distance~$u$ as thickness. The volume of
 this slice is bounded from above by the volume of the cylinder with
 base $\delta B^{d-1}$ and thickness $u$. Thus we obtain
 \begin{equation*}
  \vol(B(x,\delta)\cap B(y,\delta)) \geq \vol(\delta\ball) - u
\vol(\delta B^{d-1}) = 
\vol(\delta\ball) \lr{ 1 - u \frac{ 
\vol(\delta B^{d-1})}{\vol(\delta\ball)}}.   
 \end{equation*}
Applying Lemma~\ref{lem:bou} we obtain
$$
\frac{ \vol(\delta B^{d-1})}{\vol(\delta\ball)}=
\frac{ \vol(B^{d-1})}{\delta\vol(\ball)}
\leq \frac 1 \delta \sqrt{\frac{d+1}{2\pi}},
$$
thus  by the choice of $u\leq \sqrt{2\pi} t\delta/\sqrt{d+1} $ 
we conclude that
$$
u\frac{ \vol(\delta B^{d-1})}{\vol(\delta\ball)}
\leq  \frac{\sqrt{2\pi}t\delta
\sqrt{d+1}}{\delta\sqrt{2\pi}\sqrt {d+1}}\leq t,
$$
and the proof is complete. 
\end{proof}
We turn to the 
\begin{proof}[Proof of Lemma~\ref{lem:vempala}]
Let $x\in A_1$ and $y\in A_2$ be in $\Omega$, and suppose that their
  distance is at most $t\delta \sqrt{2 \pi/\lr{d+1}}$.
  Simple set theoretic reasoning shows that
  \begin{align*}
\vol(B(x,\delta)\cap B(y,\delta)\cap \Omega)& 
\geq \vol(B(x,\delta)\cap
\Omega) - \vol(B(x,\delta)\setminus B(y,\delta)) \\
&\geq  \vol(B(x,\delta)\cap
\Omega) - \vol(B(x,\delta)\setminus (B(x,\delta)\cap  B(y,\delta)))
\\
&= \vol(B(x,\delta)\cap \Omega) - \vol(\delta\ball) 
+ \vol(B(x,\delta)\cap  B(y,\delta)).
  \end{align*}
Since $l$ is a lower bound for the conductance $l(x)$ we have that
$$
\vol(B(x,\delta)\cap \Omega)\geq l \vol(B(x,\delta))= l
\vol(\delta\ball).
$$ 
Taking this into account and using~(\ref{eq:1}) we
end up with
\begin{align*}
 \vol(B(x,\delta)\cap B(y,\delta)\cap \Omega)& \geq l
 \vol(\delta\ball) - \vol(\delta\ball) + (1-t) \vol(\delta\ball) \\
& = (l-t)  \vol(\delta\ball). 
\end{align*}
In probabilistic terms this rewrites as
$Q_\delta(x, B(x,\delta)\cap B(y,\delta)\cap \Omega) \geq l-t$, and
similarly $Q_\delta(y, B(x,\delta)\cap B(y,\delta)\cap \Omega) \geq
l-t$.
Now, if $A\subset\Omega$ is any measurable subset with complement
$A^c$ then for $x\in A$ and $y\in A^c$ we obtain 
$$
B(x,\delta) \cap B(y,\delta)\cap\Omega \subset
\lr{B(x,\delta) \cap A^c \cap \Omega} 
\bigcup \lr{B(y,\delta) \cap A \cap \Omega} ,
$$
%E  In der letzten Formel ist ein c verschoben worden! 
which in turn yields $Q_\delta(x,A^c) + Q_\delta(y,A)\geq l-t$, but
%E auch in der letzten Formel ist ein c gewandert! 
this contradicts the definition of the sets $A_1$ and $A_2$. Hence any
two points from $A_1$ and $A_2$, respectively,  must have distance
larger than  $t\delta \sqrt{2 \pi/\lr{d+1}}$, and the proof is complete.
\end{proof}

\subsubsection*{Properties of the related Metropolis method}
\label{sec:metprop}
We analyze  Metropolis Markov chains which are based
on the ball walk, introduced above, for some appropriately chosen
$\delta$. As it will turn out, the related Metropolis chains are
\emph{perturbations} of the underlying ball walk, and its properties,
as established in Propositions~\ref{prop:ueqd} and~\ref{pro:phi}
extend in a natural way.

For $\rho \in \rad(\Omega)$ we define the \emph{acceptance
  probabilities} as
\begin{equation}
  \label{eq:alpha}
  \alph(x,y):= \min\set{1,\frac{\rho(y)}{\rho(x)}}.
\end{equation}
The corresponding Metropolis kernel is given by
\begin{equation} 
  \label{mk}
  \krd(x,dy):= 
  \alph(x,y) Q_\delta(x,dy) 
  + (1 - \int_{}\alph(x,y)Q_\delta(x,dy))\delta_x(dy).
\end{equation}
Note that for $x \notin A$ we obtain 
$$
\krd (x,A) =
\int_A \alph (x,y) \, Q_\delta (x, dy) = 
\frac 1 {\vol (\delta \ball)} \, \int_{A\cap B(x,\delta)} 
\alph (x,y) \, dy .
$$
% For the convenience of the reader
Below we sketch a single Metropolis~\KwDet 
from the present position~$x\in\Omega$ with kernel
    $\krd(x,\cdot)$. 
The procedure~{\bf Ball-walk-step} was described in
Figure~\ref{fig:bbb}.

\begin{figure}[h]
  \centering
\begin{procedure}[H]
  \caption{Metropolis-step($x,\rho,\delta$)}
\SetLine 
\Input{current position $x$, $\delta>0$, function $\rho$\;}
\Output{next position\;}
\KwProp{$y := \text{\bf Ball-walk-step}(x,\delta)$}\;
\KwAcc{}

\uIf{$\rho(y)\geq \rho(x)$}{\Return{$y$}}%{\Return{$x$}}
\uElseIf{$\rho(y) \geq {\bf rand()}\cdot \rho(x)$}{\Return{$y$}}
\Else{\Return{$x$}}
\end{procedure}  
  \caption{Schematic view of the Metropolis step. Note that the Acceptance step results in an
    acceptance probability of $\alph(x,y)=\min\set{1,\rho(y)/\rho(x)}$.}
  \label{fig:ccc}
\end{figure}
We start with the following observation.
\begin{lemma}\label{lem:beta}
Let $\alpha$ be the Lipschitz constant in $\rad(\Omega)$ and  $\beta:=
\exp(-\alpha\delta)$. 
  Uniformly for $\rho\in\rad(\Omega)$ the following bound for the
  related Metropolis chain holds true: 
  \begin{equation}
    \label{eq:alb}
\krd(x,dy) \geq \beta Q_\delta(x, dy).
  \end{equation}
\end{lemma}
\begin{proof}
Let $A\subset\Omega$. If $\dist(x,A)>\delta$  then there is nothing to
prove.
Otherwise, for $y\in A\cap B(x,\delta)$  we find
from~(\ref{eq:dens-class}) and~(\ref{eq:alpha}) that 
\begin{equation*}
 \alph(x,y)\geq \exp(-\alpha\norm{x - y}{2})
\geq e^{-\alpha\delta}=\beta.     
\end{equation*}
By definition of the transition kernel $\krd$ from~(\ref{mk}) we can
use $\beta$ to bound
$$
\krd(x,A)\geq \min\set{\alph(x,y),\ y\in A\cap B(x,\delta)}
Q_\delta(x, A) \geq \beta Q_\delta(x, A).  
$$ 
The proof is complete.
\end{proof}
The assertion of Proposition~\ref{prop:ueqd} extends to the family of
Metropolis chains as follows% , 
% quantifying similar results from~\cite{MR1399158}
.

\begin{proposition}[{cf.~\cite[Prop.~1]{MR1738303}}]\label{pro:uue1}
Let $Q_{\delta}$ be the ball walk from~(\ref{eq:qloc}) on
%  a compact set 
$\Omega$.
%  with lower bound $l$ for the local conductance of
%  $Q_{\delta/2}$. % with $\delta\leq1/\sqrt{d+1}$. 
For each $\rho\in\rad(\Omega)$ and $\delta\leq D$ the corresponding 
Metropolis chains from~(\ref{mk}) are
uniformly ergodic and reversible with respect to the related $\mur$.
\end{proposition}

\begin{proof}
Reversibility with respect to $\mur$ is clear 
by the choice of the function~$\alph$. To
prove uniform ergodicity, 
let $\beta$ be from Lemma~\ref{lem:beta} and $c$
  from~(\ref{eq:ueball}). %Set $\eta:= 1 - \beta^{n_{0}}c$.
  As established in Lemma~\ref{lem:beta} we have $\krd(x,dy)\geq
  \beta Q_{\delta}(x,dy)$. It is easy to see, and was established
  in~\cite[Proof of Thm.~2]{MR1738303}, that this extends to all
  iterates as
$$
\krd^{n}(x,dy)\geq   \beta^{n} Q^{n}_{\delta}(x,dy).
$$
Recall that under the assumptions made, 
the ball walk is uniformly ergodic, and
from Proposition~\ref{prop:ueqd} we obtain  $n_{0}$ such that for all
$x\in\Omega$ we have
\begin{equation}
  \label{eq:unifbound}
\krd^{n_{0}}(x,A)\geq  \beta^{n_{0}}c \nu(A),\quad A\subset \Omega,  
\end{equation}
proving uniform ergodicity.
\end{proof}

\begin{rem}\label{rem:unifb}
Notice that~(\ref{eq:unifbound}) is obtained with 
right hand side \emph{uniformly} for
all $\rho\in\rad(\Omega)$, a fact which will prove useful later.
\end{rem}

Finally we prove lower bounds for the conductance of the
Metropolis chains. 

\begin{theorem}\label{thm:met-cond}
Let $(\krd,\mur)$ be the Metropolis chain based
on the local ball walk $(Q_\delta,\muo)$ 
and let $\phi(\krd,\mur)$ be its conductance, where 
$\rho\in\rad(\Omega)$.
Let $l$ be a lower bound for the local conductance of $Q_{\delta}$.
For  $\rho\in\rad(\Omega)$  we have 
\begin{equation}
\label{eq:conductancelb}
\phi(\krd,\mur) \geq  
\frac{l e^{-\alpha \delta}}{8} 
\min\set{\sqrt{\frac \pi 2}\frac{l\delta}{D \sqrt{d +1}},1},
\end{equation}
where $D$ is the diameter of $\Omega$. 
\end{theorem}
\begin{rem}
  As mentioned above, Proposition~\ref{pro:phi} is a special case of
Theorem~\ref{thm:met-cond} for $\alpha=0$.
\end{rem}
The proof of Theorem~\ref{thm:met-cond} will be based on 
Lemma~\ref{lem:vempala} for the underlying
ball walk, specifying $t:= l/2$.
This extends to the Metropolis walk as follows.
\begin{lemma}\label{cor:vemp}
  Let $\alpha$
  from~(\ref{eq:dens-class}) and $l$ be the local conductance of the
  ball walk. We let $\beta:= \exp(-\alpha\delta)$.
For $A\subset \Omega$ we assign
 \begin{align}
T_1 &:= \set{x\in A,\quad  \krd(x,A^c)< \frac{\beta l}{4}}\subset
    A\\
  T_2 &:= \set{y\in A^c,\quad  \krd(y,A)< \frac{\beta l}{4}}\subset
  A^c.
  \end{align}
Then $d(T_1 ,T_2)>\delta l \sqrt{{\pi}/\lr{2d+2}}$.
 \end{lemma}
 \begin{proof}
It is enough to prove $T_1\subset A_1$ and $T_2\subset A_2$.
If $x\in T_1$ then Lemma~\ref{lem:beta} implies
$\krd(x,A^c) <\beta {l}/{4}$, hence 
$$
Q_\delta (x, A^c) \leq \frac 1 \beta \krd(x,A^c) \leq  \frac{l}{4}.
$$
The other inclusion is proved similarly.
 \end{proof}
We turn to the 
\begin{proof}[Proof of Theorem~\ref{thm:met-cond}]
Let $A\subset \Omega$ be the set for which the conductance is
attained. We assign sets $T_1$ and $T_2$ as in
Lemma~\ref{cor:vemp} and distinguish two cases. 
If $\mur (T_1)<\mur (A)/2$ \emph{or}
$\mur (T_2)<\mur (A^c)/2$, 
then the
estimate~(\ref{eq:conductancelb}) follows easily. 
For instance, if  $\mur(T_1)<\mur(A)/2$ then 
\begin{multline*}
  \int_A \krd(x,A^c)\mur(dx) \geq  \int_{A\setminus T_1}
  \krd(x,A^c)\mur(dx)\\
\geq \frac{\beta l}{4}\mur(A\setminus T_1)\geq 
\frac{ \beta l}{8}\mur(A)\geq
 \frac{\beta l}{8} \min\set{\mur(A),\mur(A^c)},
\end{multline*}
%   The choice of $t:= l/2c$ yields
% $$
%  \int_A \krd(x,A^c)\mu(dx)\geq \frac{\beta l}{8}\mu(A),
% $$
 thus $ \phi(\krd,\mur)\geq \beta l/8$ in this case, which
proves~(\ref{eq:conductancelb}).
%  under condition~(\ref{eq:condition}).

Otherwise we have $\mur(T_1)\geq \mur(A)/2$ \emph{and}
$\mur(T_2)\geq \mur(A^c)/2$. In this case we
apply an isoperimetric inequality, 
see~\cite[Thm.~4.2]{MR2178341} to the triple
$(T_1,T_2,T_3)$ with 
$T_3:= \Omega \setminus (T_1 \cup T_2)$ to conclude 
that
\begin{equation}
  \label{eq:mu123}
\mur(T_3)\geq \frac{2 d(T_1,T_2)}{D}\min\set{\mur(T_1),
\mur(T_2)},
\end{equation}
hence under the size constraints in this case it holds true that
\begin{equation}
  \label{eq:mu123f}
  \mur(T_3)\geq
  \frac{d(T_1,T_2)}{D}\min\set{\mur(A),\mur(A^c)}.
\end{equation}
Using the reversibility of the Metropolis 
chain $(\krd,\mur)$ we have
$$
\int_A \krd(x,A^c)\mur(dx)= \int_{A^c} \krd(y,A)\mur(dy),
$$
which implies
\begin{align*}
\int_A \krd(x,A^c)\mur(dx)&= \frac 1 2 \lr{\int_A \krd(x,
  A^c )\mur(dx)+ \int_{A^c} \krd(y,A)\mur(dy) }  \\
& \ge  \frac 1 2 \lr{ \int_{A\cap T_3} \krd(x,
  A^c )\mur(dx)+ \int_{A^c \cap T_3} \krd(y,A)\mur(dy) }\\
&\geq \frac 1 2 \lr{ \frac{\beta l }{4} \mur(A \cap T_3) +
  \frac{\beta l}{4} \mur(A^c  \cap T_3) }\\
&= \frac{\beta l }{8}\lr{\mur(A \cap T_3) 
+\mur(A^c  \cap T_3) }=
\frac{\beta l}{8} \mur(T_3).
\end{align*}
Since by Lemma~\ref{cor:vemp} we can bound
$d(T_1,T_2)\geq \delta l \sqrt{{\pi}/\lr{2d+2}}$ 
we use~(\ref{eq:mu123f}) to
complete the proof.
\end{proof}

If we restrict ourselves to Metropolis chains on $\ball$, then
Lemma~\ref{lem:l-bound}  provides a lower bound for 
the local conductance which is independent of the
dimension~$d$. 
As a simple consequence of Theorem~\ref{thm:met-cond} we 
then obtain the following
\begin{corollary}
\label{cor2} 
Assume that $\rho \in \rad(\ball)$ and 
$\delta \le (d+1)^{-1/2}$. 
Then we obtain 
$$
\phi(K_{\rho, \delta} , \mur) \ge\sqrt{\frac \pi 2}
\frac{9 \delta}{1600 \sqrt{d + 1}} e^{-\alpha \delta} .
$$
To maximize $\phi$ we define
\begin{math} % \label{eq:max} 
\delta^* = \min\set{{1}/{\sqrt{d+1}},1 /\a }
\end{math} 
and obtain 
$$
\phi(K_{\rho, \delta^*}, \mur)   
\ge 0.0025 \, \frac{1}{\sqrt{d+1}} 
\min\set{\frac{1}{\sqrt{d+1}},\frac 1 \a } .
$$
\end{corollary}

\subsubsection*{Error bounds}
\label{sec:er-mc}

For the class $\fad(\Omega)$ 
the above lower conductance bound~(\ref{eq:conductancelb}) 
will yield an error estimate for the problem~(\ref{eq02}).

Let $S_n^\delta$ be the 
estimator based on a sample of the local Metropolis Markov
chain with transition $K_{\rho,\delta}$, starting at zero.
To estimate its error
we combine the estimates of the conductance of $K_{\rho,\delta}$
with two results, partially known from the literature.
To formulate the results we note the following. 
The Markov kernel 
$K_{\rho, \delta}$ is reversible with respect to 
$\mur$ and hence induces a self-adjoint operator
$$
K_{\rho, \delta} : L_2 (\Omega,\mur) \to L_2 (\Omega,\mur) .
$$
The spectrum $\sigma (K_{\rho,\delta})$ is contained 
in $[-1, 1]$ and $1 \in \sigma(K_{\rho, \delta})$ 
and we are interested in the second largest eigenvalue 
$$
\beta_{\rho, \delta} := \sup \{ \sigma \in \sigma 
(K_{\rho, \delta}) \mid \sigma \not= 1 \}
$$
 of $K_{\rho,\delta}$. This is motivated by the  extension of  a result 
from~\cite[Cor.~1]{MR1738303} about the worst case
error of $S_n^\delta$, uniformly for $(f,\rho)\in\fad(\Omega)$.  
\begin{lemma}
\label{le:mathescharf} 
$$
\lim_{n \to \infty } \sup_{(f,\rho)\in\fad(\Omega)} 
e(S_n^\delta, (f,\rho))^2 \cdot n =
\sup_{\rho\in\rad(\Omega)}\frac{1+ \beta_{\rho, 
\delta}}{1- \beta_{\rho, \delta}} . 
$$
\end{lemma} 
The proof is given in the appendix. 
For Markov chains which start according to the invariant distribution
$\mur$ the bound is similar, but  more explicit and was given
in~\cite{SOK} and~\cite[Thm.~1.9]{MR1238906}.

The relation of the second largest
eigenvalue~$\beta_{\rho, \delta}$ to the conductance 
is given in 

\begin{lemma}[Cheeger's Inequality, 
see~\cite{MR1025467,MR930082,MR1238906}]
\label{le:cheeger}
$$
\lambda_{\rho,\delta} :=  1 - \beta_{\rho, \delta} \geq 
\phi^{2}(K_{\rho, \delta}, \mur)/2.
$$
\end{lemma}

We are ready to state %  and prove
our main result for 
the Metropolis algorithm $S_n^\delta$, based on the Markov chain 
$K_{\rho, \delta}$, for the class 
$\fad (\ball)$, i.e., when $\Omega\subset \R^{d}$ is the Euclidean unit ball. 
\begin{theorem} 
\label{th5} 
Let $S_n^\delta=\frac 1 n \sum_{j=1}^{n}f(X_{j})$ be the 
estimator based on a sample~$(X_{1},\dots,X_{n})$ of the local Metropolis Markov
chain with transition $K_{\rho, \delta}$, 
where $\delta \le (d+1)^{-1/2}$.
Then
\begin{equation}
  \label{eq:th5}
 \lim_{n \to \infty} \sup_{(f,\rho)\in\fad(\ball)}  
e(S_n^\delta, (f,\rho) ) ^2 \cdot n 
\le \frac{8\cdot 1600^{2}}{81\pi}(d +1)\cdot \frac{e^{2 \alpha \delta}}
{\delta^{2}} . 
\end{equation}
 Again we may choose 
$
\delta^* = \min\set{(d+1)^{-1/2},\alpha^{-1}}
$
and obtain 
\begin{equation} 
\label{tract} 
\lim_{n \to \infty} \sup_{(f,\rho)\in\fad(\ball)} 
e(S_n^{\delta^*} , (f,\rho) ) ^2 \cdot n 
\le 594700 \cdot (d+1)\max\set{d+1,{\alpha^{2}}}. 
\end{equation} 
\end{theorem} 
\begin{proof} 
This follows from Corollary~\ref{cor2}, and 
 Lemmas~\ref{le:mathescharf} and~\ref{le:cheeger}. 
\end{proof} 

\section{Summary}
\label{sec:sum}

Let us discuss our findings.  %   in some detail.
The results from Section~\ref{s2} clearly indicate that the
superiority of Metropolis algorithms upon 
simpler (non-adaptive) Monte Carlo methods
does not hold in general. Specifically, it does not hold 
for the large classes $\fco$ of input without
additional structure.

On the other hand, for the class~$\fad(\ball)$, specific Metropolis
algorithms that are based on local underlying walks are superior to
all non-adaptive methods. 
Even more,   %  , as formula~\eqref{tract} indicates,
on~$\ball$  %   the problem is \emph{tractable} 
%   in $d$ and $\alpha$:
the cost of the algorithm~$S_n^{\delta^*}$, roughly 
given by the number $n$ of evaluations of $\rho$ and $f$, 
increases like a polynomial in $d$ and $\alpha$.
More 
precisely, according to~\eqref{tract}, the asymptotic constant 
$\lim_{n \to \infty} e(S_n^{\delta^*} , \fad(\ball) ) ^2~\cdot~n$
is bounded by a constant times~$\max\set{d^{2}, d\alpha^{2}}$, 
i.e., the complexity grows polynomially in $d$ and $\alpha$
and, for fixed $d$, increases (at most) as $\alpha^{2}$. 
If we only allow non-adaptive methods then this asymptotic constant,
again for fixed $d$,  increases at least as $\alpha^{d}$,
see~\eqref{lo9}. 

%E  Hier Modifikationen. 
We believe that this problem is \emph{tractable} in the sense that 
the number of function values to achieve an error $\e$ can be bounded
by 
\begin{equation}   \label{tract2} 
n(\e , \fad(\ball)   ) \,  \le \,  C \,  \e^{-2} \,  d \,  \max ( d, \alpha^2) .
\end{equation} 
We did not prove \eqref{tract2}, however, since Theorem 5 is only a
statement for large $n$. 

Notice
that according to Theorem~\ref{th5} the size~$\delta^{\ast}$ of the
underlying balls walk needs to be adjusted both to 
the spatial dimension~$d$ and the
Lipschitz constant~$\alpha$.

The analysis of the Metropolis algorithm is based on properties of the
underlying ball walk; in particular we establish uniform ergodicity of
the ball walk for convex bodies~$\Omega\subset \R^{d}$. Also, based
on conductance arguments,  we provide lower bounds for the spectral gap
of the ball walk.   % : If $\delta\sim 1/\sqrt{d}$, then
% Proposition~\ref{pro:phi} together 
% with Lemma~\ref{le:cheeger} show that this is independent
% of the dimension

As a consequence,  in the case~$\alpha=0$ the estimate~(\ref{eq:th5})  provides an
error bound for the ball 
walk $(Q_{\delta},\mu)$, which is asymptotically of the form 
$ e(S_{n}^{\delta},L_{2}(\ball,\mu))\leq C \delta^{-1}
(d/n)^{1/2}$.%  This complements the heuristic considerations
% from~\cite[Example~1]{MR1738303}.

The results  extend in a similar way to any family 
$\Omega_d \subset \R^{d}$ for which
the underlying local ball walk $Q_{\delta}$ has 
(for $\delta \le \delta_d$) 
a non-trivial lower bound for the
local conductance that is independent of the dimension.

Finally, from the results of Section~\ref{s2} we can conclude that adaption
does not help much for the classes $\fco$. 
Hence we have new results concerning the \emph{power of
adaption}, see~\cite{MR1408328} for a survey of earlier results, in
particular that it may help to break the \emph{curse of
dimensionality} for the classes $\fad(\ball)$. 

\appendix
\section{Proof of Lemma~\ref{le:mathescharf}}
% \label{app}

Lemma~\ref{le:mathescharf} extends the bound 
from~\cite[Thm.~1]{MR1738303}, which deals with a single uniformly
ergodic chain. It was obtained from  on a contraction
property, as stated in~\cite[Prop.~1]{MR1738303}.
The goal of the
present analysis is to establish this asymptotic result 
\emph{uniformly} for all Metropolis chains with density from
$\rad(\Omega)$, by showing that this contractivity holds true uniformly% , which in turn allows to
% extend the proof of Theorem~1 in~\cite{MR1738303}
.

\subsection*{Contractivity of the Markov operator}
We assign to each transition kernel $K$ on $\Omega$ with corresponding invariant
distribution $\mu$ the bounded linear mapping $P$, given by
\begin{equation}
  \label{eq:prd}
(P f)(x) := \int f(y) K(x,dy).
\end{equation}
Also we let $E$ denote the mapping which assigns any integrable
function its expectation as a constant function
$
E(f)\colon= \int_\Omega f(x) \, \mu(dx).
$
{F}or each $K$ the mapping $P - E$ is bounded in
$L_{\infty}(\Omega,\mu)$, with norm less than or  equal to one and we
shall strengthen this uniformly for kernels $\krd$ with
$\rho\in\rad(\Omega)$.
Within this operator context 
\emph{uniform ergodicity} is equivalent to a specific
form of quasi-compactness, namely there are $0<\eta<1$ and $n_{0}\in\N$
for which
\begin{equation}\label{eq-infty-con}
\norm{P^{n} - E\colon L_{\infty}(\Omega)
\to L_{\infty}(\Omega)}{}\leq\eta,\ \text{for $n\geq
n_{0}$.}
\end{equation}
We first show that reversibility allows to transfer this to
the spaces~$L_{1}(\Omega,\mu_{\rho})$.
\begin{lemma} \label{lem:infty1}
Suppose that the transition kernel $K$ with corresponding
  mapping $P$  is reversible. Then for all $n\in\N$ we have
  \begin{equation}
    \label{eq:1infty}
    \norm{P^{n} - E\colon L_{1}(\Omega,\mu)
    \to L_{1}(\Omega,\mu)}{}
    \leq  \norm{P^{n} - E\colon L_{\infty}(\Omega,\mu)
    \to L_{\infty}(\Omega,\mu)}{}.
  \end{equation}
% Consequently, if $K$ is uniformly ergodic and reversible, then there
% are $n_{0}\in\N$ and $\eta<1$ such that 
% \begin{equation}
%   \label{eq:1eta}
% \norm{P^{n_{0}} - E\colon L_{1}(\ball,\mu)
% \to L_{1}(\ball,\mu)}{}\leq \eta.
% \end{equation}
  \end{lemma}
  \begin{proof}
    If $K$ is reversible, then so are all iterates $K^{n}$. Thus for
    arbitrary functions $f\in L_{1}(\Omega,\mu)$ and $h\in
    L_{\infty}(\Omega,\mu)$ we have, using the scalar product on
    $L_{2}(\Omega,\mu)$, that
$$
\scalar{(P^{n}- E)f}{h}= \scalar{f}{(P^{n}- E)h}.
$$
Consequently, for any $f\in L_{1}(\Omega,\mu)$ we have
\begin{align*}
  \norm{(P^{n} - E) f}{1} &= \sup_{\norm{h}{\infty}\leq 1}
  \abs{\scalar{(P^{n}- E)f}{h}} =
  \sup_{\norm{h}{\infty}\leq 1} \abs{\scalar{f}{(P^{n}- E)h}}  \\
  &\leq \norm{f}{1} \sup_{\norm{h}{\infty}\leq 1} \norm{(P^{n}-
    E)h}{\infty},
\end{align*}
from which the proof can be completed.
\end{proof}

\begin{proposition}\label{pro:unbound}
% Let $\Omega\subset\R^{d}$ be a compact set, 
% and suppose that there is  a lower bound $l>0$ 
% for the local conductance of $Q_{\delta/2}$.
For any convex body $\Omega \subset \R^d$ 
there are an
integer $n_{0}$ and a constant $0<\eta<1$ such that uniformly for
$\rho\in\rad(\Omega)$ we have
\begin{equation}
  \label{eq:uue1}
  \norm{\prd^{n_{0}} - E\colon L_{1}(\Omega,\mu_{\rho})\to
    L_{1}(\Omega,\mu_{\rho})}{}\leq \eta.
\end{equation}
\end{proposition}

\begin{proof}
This is  an immediate consequence of the
bound~(\ref{eq:unifbound}). As mentioned in Remark~\ref{rem:unifb}
uniform ergodicity was established uniformly for $\rho\in\rad(\Omega)$.
 It is well known (see~\cite[Thm.~16.2.4]{Meyn-book}) that this
implies % the assertion~(\ref{eq-infty-con}).
% because this yields
that there is an $\eta<1$ such that uniformly for
$\rho\in\rad(\Omega)$ we have
\begin{equation}\label{eq-infty-con2}
\norm{\prd^{n_{0}} - E\colon L_{\infty}(\Omega)
\to L_{\infty}(\Omega)}{}\leq\eta,\ \text{for $n\geq
n_{0}$.}
\end{equation}
In the light of Lemma~\ref{lem:infty1} this yields~(\ref{eq:uue1}).
\end{proof}

Finally we sketch the 
\begin{proof}[Proof of Lemma~\ref{le:mathescharf}]
Using Proposition~\ref{pro:unbound} we can extend the
proof of~\cite[Thm.~1]{MR1738303}. 
In particular, the bounds from Eq.~(13)--(15)
in~\cite{MR1738303} tend to zero uniformly for
$\rho\in\rad(\Omega)$. Moreover, starting at zero, 
after one step according to the underlying ball walk, the (new)
initial distribution is uniformly bounded  with respect to the uniform
distribution on $\Omega$, hence also with respect to~$\mur$,
such that we establish the asymptotics in Lemma~\ref{le:mathescharf}.  
\end{proof}

\medskip
\noindent
{\bf Acknowledgment:} \
We thank two anonymous referees and Daniel Rudolf for their comments. 
%E  Ackno ist neu. 

%\cite{MR2260070,MR2172842}

\bibliographystyle{plain}
%\bibliography{ref,mybib}

\begin{thebibliography}{10}

\bibitem{MR2260070}
Christophe Andrieu and {\'E}ric Moulines.
\newblock On the ergodicity properties of some adaptive {MCMC} algorithms.
\newblock {\em Ann. Appl. Probab.}, 16(3):1462--1505, 2006.

\bibitem{103439}
David Applegate and Ravi Kannan.
\newblock Sampling and integration of near log-concave functions.
\newblock In {\em STOC '91: Proceedings of the twenty-third annual ACM
  symposium on Theory of computing}, pages 156--163, New York, NY, USA, 1991.
  ACM Press.

\bibitem{MR2172842}
Yves~F. Atchad{\'e} and Jeffrey~S. Rosenthal.
\newblock On adaptive {M}arkov chain {M}onte {C}arlo algorithms.
\newblock {\em Bernoulli}, 11(5):815--828, 2005.

\bibitem{Bachvalov}
N.~S. Bahvalov.
\newblock Approximate computation of multiple integrals.
\newblock {\em Vestnik Moskov. Univ. Ser. Mat. Meh. Astr. Fiz. Him.},
  1959(4):3--18, 1959.

\bibitem{B/D06}
F.~Bassetti and P.~Diaconis.
\newblock Examples comparing importance sampling and the {M}etropolis
  algorithm.
\newblock {\em to appear Illinois J. of Math.}, 2006.

\bibitem{10.1109/5992.814660}
Isabel Beichl and Francis Sullivan.
\newblock The {M}etropolis algorithm.
\newblock {\em Computing in Science and Engineering}, 2(1):65--69, 2000.

\bibitem{10.1109/MCSE.2006.27}
Isabel Beichl and Francis Sullivan.
\newblock Guest editors' introduction: Monte {C}arlo methods.
\newblock {\em Computing in Science and Engineering}, 8(2):7--8, 2006.

\bibitem{MR2013000}
Nicolas Bourbaki.
\newblock {\em Functions of a real variable}.
\newblock Elements of Mathematics (Berlin). Springer-Verlag, Berlin, 2004.

\bibitem{Burenkov} 
Victor I. Burenkov.
{\em Sobolev Spaces on Domains.} 
Teubner-Texte zur Mathematik 137.
Teubner Verlag Stuttgart, 1998. 

\bibitem{MR1284987}
Alan Frieze, Ravi Kannan, and Nick Polson.
\newblock Sampling from log-concave distributions.
\newblock {\em Ann. Appl. Probab.}, 4(3):812--837, 1994.

\bibitem{Hlawka}
E.~Hlawka.
\newblock Ausf\"ullung und \"Uberdeckung konvexer K\"orper durch 
konvexe K\"orper.
\newblock {\em Mh. Math. Phys.}, 53:81--131, 1949.

\bibitem{MR1025467}
Mark Jerrum and Alistair Sinclair.
\newblock Approximating the permanent.
\newblock {\em SIAM J. Comput.}, 18(6):1149--1178, 1989.

\bibitem{MR1318794}
R.~Kannan, L.~Lov{\'a}sz, and M.~Simonovits.
\newblock Isoperimetric problems for convex bodies and a localization lemma.
\newblock {\em Discrete Comput. Geom.}, 13(3-4):541--559, 1995.

\bibitem{MR797411}
Ulrich Krengel.
\newblock {\em Ergodic theorems}, volume~6 of {\em de Gruyter Studies in
  Mathematics}.
\newblock Walter de Gruyter \& Co., Berlin, 1985.

\bibitem{MR930082}
Gregory~F. Lawler and Alan~D. Sokal.
\newblock Bounds on the {$L\sp 2$} spectrum for {M}arkov chains and {M}arkov
  processes: a generalization of {C}heeger's inequality.
\newblock {\em Trans. Amer. Math. Soc.}, 309(2):557--580, 1988.

\bibitem{MR1238906}
L.~Lov{\'a}sz and M.~Simonovits.
\newblock Random walks in a convex body and an improved volume algorithm.
\newblock {\em Random Structures Algorithms}, 4(4):359--412, 1993.

\bibitem{olm}
Peter Math{\'e}.
\newblock The optimal error of {M}onte {C}arlo integration.
\newblock {\em J. Complexity}, 11(4):394--415, 1995.

\bibitem{MR1738303}
Peter Math{\'e}.
\newblock Numerical integration using {M}arkov chains.
\newblock {\em Monte Carlo Methods Appl.}, 5(4):325--343, 1999.

\bibitem{Meyn-book}
S.~P. Meyn and R.~L. Tweedie.
\newblock {\em Markov chains and stochastic stability}.
\newblock Springer-Verlag London Ltd., London, 1993.

\bibitem{NOV}
Erich Novak.
\newblock {\em Deterministic and stochastic error bounds in numerical
  analysis}.
\newblock Lect. Notes Math. 1349. Springer-Verlag, Berlin, 1988.

\bibitem{MR1319050}
Erich Novak.
\newblock The real number model in numerical analysis.
\newblock {\em J. Complexity}, 11(1):57--73, 1995.

\bibitem{MR1408328}
Erich Novak.
\newblock On the power of adaption.
\newblock {\em J. Complexity}, 12(3):199--237, 1996.

\bibitem{10.1109/MCSE.2006.30}
Dana Randall.
\newblock Rapidly mixing {M}arkov chains with applications in computer science
  and physics.
\newblock {\em Computing in Science and Engineering}, 8(2):30--41, 2006.

\bibitem{MR1399158}
G.~O. Roberts and R.~L. Tweedie.
\newblock Geometric convergence and central limit theorems for multidimensional
  {H}astings and {M}etropolis algorithms.
\newblock {\em Biometrika}, 83(1):95--110, 1996.

\bibitem{MR0172183}
C.~A. Rogers.
\newblock {\em Packing and covering}.
\newblock Cambridge Tracts in Mathematics and Mathematical Physics, No. 54.
  Cambridge University Press, New York, 1964.

\bibitem{SOK}
A.~Sokal.
\newblock Monte {C}arlo methods in statistical mechanics: foundations and new
  algorithms.
\newblock In {\em Functional integration (Carg\`ese, 1996)}, pages 131--192.
  Plenum, New York, 1997.

\bibitem{IBC}
J.~F. Traub, G.~W. Wasilkowski, and H.~Wo{\'z}niakowski.
\newblock {\em Information-based complexity}.
\newblock Academic Press Inc., Boston, MA, 1988.
\newblock With contributions by A. G. Werschulz and T. Boult.

\bibitem{vempala-lesson}
Santosh Vempala.
\newblock Lect.~17, {R}andom walks and polynomial time algorithms.
\newblock {http://www-math.mit.edu/\~{}vempala/random/course.html}, 2002.

\bibitem{MR2178341}
Santosh Vempala.
\newblock Geometric random walks: a survey.
\newblock In {\em Combinatorial and computational geometry}, volume~52 of {\em
  Math. Sci. Res. Inst. Publ.}, pages 577--616. Cambridge Univ. Press,
  Cambridge, 2005.

\end{thebibliography}

\def\cprime{$'$} \def\cprime{$'$}

\end{document}